\magnification=\magstep1
\newif\ifsect\newif\iffinal
\secttrue\finaltrue
\def\strutdepth{\dp\strutbox}

\def\lsimb#1{\vadjust{\vtop to0pt{\baselineskip\strutdepth\vss
	\llap{\ttt\string #1\ }\null}}}
\def\rsimb#1{\vadjust{\vtop to0pt{\baselineskip\strutdepth\vss
	\line{\kern\hsize\rlap{\ttt\ \string #1}}\null}}}
\def\ssect #1. {\bigbreak\indent{\bf #1.}\enspace\message{#1}}
\def\smallsect #1. #2\par{\bigbreak\noindent{\bf #1.}\enspace{\bf #2}\par
	\global\parano=#1\global\eqnumbo=1\global\thmno=1
	\nobreak\smallskip\nobreak\noindent\message{#2}}
\def\thm #1: #2{\medbreak\noindent{\bf #1:}\if(#2\thmp\else\thmn#2\fi}
\def\thmp #1) { (#1)\thmn{}}
\def\thmn#1#2\par{\enspace{\sl #1#2}\par
        \ifdim\lastskip<\medskipamount \removelastskip\penalty 55\medskip\fi}

\def\qedn{\thinspace\null\nobreak\hfill\hbox{\vbox{\kern-.2pt\hrule height.2pt 
depth.2pt\kern-.2pt\kern-.2pt \hbox to2.5mm{\kern-.2pt\vrule width.4pt
\kern-.2pt\raise2.5mm\vbox to.2pt{}\lower0pt\vtop to.2pt{}\hfil\kern-.2pt
\vrule width.4pt\kern-.2pt}\kern-.2pt\kern-.2pt\hrule height.2pt depth.2pt
\kern-.2pt}}\par\medbreak}
\def\pf{\ifdim\lastskip<\smallskipamount \removelastskip\smallskip\fi
        \noindent{\sl Proof\/}:\enspace}
\def\itm#1{\par\indent\llap{\rm #1\enspace}\ignorespaces}
\let\smbar=\bar
\def\bar#1{\overline{#1}}
\def\forclose#1{\hfil\llap{$#1$}\hfilneg}
\def\newforclose#1{
	\ifsect\xdef #1{(\number\parano.\number\eqnumbo)}\else
	\xdef #1{(\number\eqnumbo)}\fi
	\hfil\llap{$#1$}\hfilneg
	\global \advance \eqnumbo by 1
	\iffinal\else\rsimb#1\fi}
\def\forevery#1#2$${\displaylines{\let\eqno=\forclose
        \let\neweq=\newforclose\hfilneg\rlap{$\qqquad\forall#1$}\hfil#2\cr}$$}
\newcount\parano
\newcount\eqnumbo
\newcount\thmno
\newcount\versiono
\versiono=0
\def\neweqt#1$${\xdef #1{(\number\parano.\number\eqnumbo)}
	\eqno #1$$
	\iffinal\else\rsimb#1\fi
	\global \advance \eqnumbo by 1}
\def\newthmt#1 #2: #3{\xdef #2{\number\parano.\number\thmno}
	\global \advance \thmno by 1
	\medbreak\noindent
	\iffinal\else\lsimb#2\fi
	{\bf #1 #2:}\if(#3\thmp\else\thmn#3\fi}
\def\neweqf#1$${\xdef #1{(\number\eqnumbo)}
	\eqno #1$$
	\iffinal\else\rlap{$\smash{\hbox{\hfilneg\string#1\hfilneg}}$}\fi
	\global \advance \eqnumbo by 1}
\def\newthmf#1 #2: #3{\xdef #2{\number\thmno}
	\global \advance \thmno by 1
	\medbreak\noindent
	\iffinal\else\llap{$\smash{\hbox{\hfilneg\string#1\hfilneg}}$}\fi
	{\bf #1 #2:}\if(#3\thmp\else\thmn#3\fi}
\def\bititolo{\empty}
\gdef\begin #1 #2\par{\xdef\titolo{#2}
\ifsect\let\neweq=\neweqt\else\let\neweq=\neweqf\fi
\ifsect\let\newthm=\newthmt\else\let\newthm=\newthmf\fi
\centerline{\titlefont\titolo}
\if\bititolo\empty\else\medskip\centerline{\titlefont\bititolo}
\xdef\titolo{\titolo\ \bititolo}\fi
\bigskip
\centerline{\bigfont by \autore}
\if\istituto!\else\bigskip
\centerline{\istituto}
\centerline{\indirizzo}\fi
\medskip
\centerline{#1~\anno}
\bigskip\bigskip
\ifsect\else\global\thmno=1\global\eqnumbo=1\fi}
\font\titlefont=cmssbx10 scaled \magstep1
\font\bigfont=cmr10 scaled \magstep1

\font\ttt=cmtt10 at 10truept
\font\eightrm=cmr8

\font\bbr=msbm10
\font\sbbr=msbm7 
\def\ca #1{{\cal #1}}
\nopagenumbers
\binoppenalty=10000
\relpenalty=10000
\let\de=\partial
\def\eps{\varepsilon}

\def\phe{\varphi}

\def\Re{\mathop{\rm Re}\nolimits}

\def\cancel#1#2{\ooalign{$\hfil#1/\hfil$\crcr$#1#2$}}
\def\void{\mathord{\mathpalette\cancel{\mathrel{\hbox{\kern0pt\raise0.8pt\hbox
	{$\scriptstyle\bigcirc$}}}}}}

\def\R{{\mathchoice{\hbox{\bbr R}}{\hbox{\bbr R}}{\hbox{\sbbr R}}
{\hbox{\sbbr R}}}}
\def\C{{\mathchoice{\hbox{\bbr C}}{\hbox{\bbr C}}{\hbox{\sbbr C}}
{\hbox{\sbbr C}}}}

\def\qqquad{\quad\qquad}

\newcount\notitle
\notitle=1
\headline={\ifodd\pageno\rhead\else\lhead\fi}
\def\rhead{\ifnum\pageno=\notitle\iffinal\hfill\else\hfill\tt Version 
\the\versiono; \the\day/\the\month/\the\year\fi\else\hfill\eightrm\titolo\hfill
\folio\fi}
\def\lhead{\ifnum\pageno=\notitle\hfill\else\eightrm\folio\hfill\autore\hfill
\fi}
\def\autore{Marco Abate}
\output={\plainoutput}
\newbox\bibliobox
\def\setref #1{\setbox\bibliobox=\hbox{[#1]\enspace}
	\parindent=\wd\bibliobox}
\def\biblap#1{\noindent\hang\rlap{[#1]\enspace}\indent\ignorespaces}
\def\art#1 #2: #3! #4! #5 #6 #7-#8 \par{\biblap{#1}#2: {\sl #3\/}.
	#4 {\bf #5}~(#6)\if.#7\else, \hbox{#7--#8}\fi.\par\smallskip}
\def\book#1 #2: #3! #4 \par{\biblap{#1}#2: {\bf #3.} #4.\par\smallskip}
\versiono=2
\overfullrule=0pt
\def\nablab{{\nabla\!}}
\def\llangle{{\langle\!\langle}}
\def\rrangle{{\rangle\!\rangle}}
\def\dmatrix{\def\normalbaselines{\baselineskip20pt\lineskip3pt
	\lineskiplimit3pt}\matrix}
\def\mapright#1{\smash{\mathop{\longrightarrow}\limits^{#1}}}
\def\mapdownr#1{\Big\downarrow\rlap{$\vcenter{\hbox{$\scriptstyle#1$}}$}}
\def\mapdownl#1{\llap{$\vcenter{\hbox{$\scriptstyle#1$}}$}\Big\downarrow}
\def\mapse#1{\llap{$\vcenter{\hbox{$\scriptstyle#1$}}$}\searrow}
\def\Wedge^#1{{\textstyle\bigwedge^{\!#1}}}

\font\bbr=msbm10
\font\sbbr=msbm7 
\def\ca #1{{\cal #1}}
\let\smb=\smbar

\def\R{{\mathchoice{\hbox{\bbr R}}{\hbox{\bbr R}}{\hbox{\sbbr R}}
{\hbox{\sbbr R}}}}
\def\C{{\mathchoice{\hbox{\bbr C}}{\hbox{\bbr C}}{\hbox{\sbbr C}}
{\hbox{\sbbr C}}}}

\def\tildeM{\tilde M} 
\def\autore{Marco Abate and Giorgio Patrizio}
\def\istituto{Seconda Universit\`a di Roma}\def\indirizzo{00133 Roma, Italy}
\def\anno{1993}
%
%
%
%
%
%
\begin September Complex Finsler metrics
 
\centerline{\it Dedicated to Prof. S. Kobayashi}
 
\smallsect 0. Introduction
 
A complex Finsler metric is an upper semicontinuous function $F\colon T^{1,0}M
\to\R^+$ defined on the holomorphic tangent bundle of a complex Finsler 
manifold~$M$, with the property that
$F(p;\zeta v)=|\zeta|F(p;v)$ for any $(p;v)\in T^{1,0}M$ and $\zeta\in\C$. 
 
Complex Finsler metrics do occur naturally in function theory of 
several variables. The Kobayashi metric introduced in 1967 ([K1]) and its 
companion the Carath\'eodory metric are remarkable examples which have become 
standard tools for anybody working in complex analysis; we refer the reader 
to~[K2,~4], [L], [A] and~[JP] to get an idea of the amazing developments in 
this area achieved in the past 25~years.
 
In general, the Kobayashi metric is not at all regular; it may even not be 
continuous. But in 1981 Lempert~[Le] proved that the Kobayashi metric of a 
bounded strongly convex domain $D$ in $\C^n$ is smooth (outside the zero 
section of $T^{1,0}D$), thus allowing in principle the use of differential 
geometric techniques in the study of function theory over strongly convex 
domains (see also Pang~[P2] for other examples of domains with smooth 
Kobayashi metric).
 
We started dealing with this kind of problems in [AP1]. In particular, [AP2] 
was devoted to the search of differential geometric conditions ensuring the 
existence in a complex Finsler manifold of a foliation in holomorphic 
disks like the one found by Lempert in strongly convex domains, where the 
disks were isometric embeddings of the unit disk $\Delta\subset\C$ endowed 
with the Poincar\'e metric. And indeed (see also~[AP3]) we found necessary and 
sufficient conditions (see also Pang~[P1] for closely related results).
In that case, because the nature of the problem required the solution of 
certain P.D.E.'s, the conditions were mainly expressed in local coordinates
somewhat hiding their geometric meaning.
 
The aim of this paper is to present an introduction to complex Finsler 
geometry in a way suitable to deal with global questions. Roughly speaking, 
the idea is to isometrically embed a complex Finsler manifold into a hermitian 
vector bundle, and then apply standard hermitian differential geometry 
techniques, in the spirit of~[K3]. Here we provide just a 
coarse outline of the procedure. Let $\tildeM$ be the complement of the zero 
section in~$T^{1,0}M$. We assume that the complex Finsler metric~$F$ 
is smooth on~$\tildeM$, and that $F$ is {\sl 
strongly pseudoconvex,} that is that the Levi form of $G=F^2$ is positive 
definite. Now let $\ca V\subset T^{1,0}\tildeM$ be the {\sl vertical 
bundle,} that is the kernel of the differential of the canonical projection
$\pi\colon T^{1,0}M\to M$. Using the Levi form of~$G$, it is easy to define a 
hermitian metric on~$\ca V$; moreover, there exists a canonical 
section~$\iota$ of~$\ca V$ giving an isometric embedding of~$\tildeM$ 
into~$\ca V$ --- that is for any~$v\in\tildeM$ the norm of~$\iota(v)$ with 
respect to the given hermitian metric on~$\ca V$ is equal to~$F(v)$.
Let $D$ be the Chern connection on~$\ca V$ associated to the metric, and 
denote by~$\ca H$ the kernel of the bundle map $X\mapsto\nablab_X\iota$. Then 
it turns out that $\ca H$ is a {\sl horizontal bundle,} that is $T^{1,0}M=\ca 
H\oplus\ca V$; furthermore, there is a canonically defined global bundle 
isomorphism $\Theta\colon\ca V\to\ca H$. Using~$\Theta$, we can transfer both 
the metric and the connection on~$\ca H$, obtaining a canonical hermitian 
structure on~$T^{1,0}\tildeM$, and the associated Chern connection preserves 
the splitting. Finally, the {\sl horizontal radial vector field} $\chi=\Theta
\circ\iota$ is a canonical isometric embedding of~$\tildeM$ into~$\ca H$. Then 
our idea is that the complex Finsler geometry of~$M$ should be described by 
using the differential geometry of the Chern connection~$D$ restricted to~$\ca 
H$, using~$\chi$ as a means of transfering informations from the tangent 
bundle to the horizontal bundle and back. For instance, the K\"ahler condition 
introduced in~[AP2] becomes the vanishing of a suitable contraction of the 
horizontal part of the torsion of~$D$ (here we say that the metric is {\sl 
weakly K\"ahler\/}); and the necessary and sufficient condition for the 
existence of complex geodesic curves (see~[AP2,~3]) are expressed by constant 
holomorphic curvature and a symmetry property of the horizontal part of the 
curvature of~$D$; cf.~Lemma~8.3.
 
This approach is in the spirit of the one developed by E.~Cartan~[C] for real 
Finsler metrics; see~[Ru1], [M], [Ch], [BC], [Be] and the forthcoming 
monograph~[AP4] for an account in modern language. On the other hand, to our 
surprise we 
were unable to find in the literature a comparable approach in the complex 
case. Rund, in~[Ru2], described the Chern connection on the horizontal bundle, 
but only in local coordinates. Fukui in~[Fu] studied the Cartan connection on 
a complex Finsler manifold, which is in general different from the Chern
connection (see~[AP4] for a comparison). Faran~[F] studied the local 
equivalence problem, without dealing with global questions.
Only Kobayashi~[K3] explicitely used the Chern connection, but he seemed 
unaware of the relevance of the horizontal component. It should be mentioned 
that we choose to work on~$\tildeM$ instead of the projectivized tangent bundle 
mainly for keeping more transparent the relationships between global objects 
and local computations (which are often simplified by consistently using the 
homogeneity of the function~$G$ and its derivatives). However the two approach 
are completely equivalent. In fact, the role of the canonical sections~$\iota$ 
and~$\chi$ in our context is analogous to the role of the tautological line 
bundle in~[K3].
We hope that our work will clarify the subject of complex Finsler geometry,
opening the way to new research in the field. 
 
The content of this paper is the
following. In sections~1 and~2 we describe in detail the construction outlined 
above of the Chern-Finsler connection.
In sections~3 and~4 we define the (2,0)-torsion, the (1,1)-torsion, the
curvature of the Chern-Finsler connection on the horizontal bundle, we derive 
the Bianchi identities and we discuss K\"ahler Finsler metrics. 
In section~5 we introduce the notion of holomorphic curvature. 
 
In sections~6 and~7  we derive the first and second variation formulas for a 
strongly pseudoconvex K\"ahler Finsler metric, giving a good example of global 
computations made using the tools introduced before. As a corollary, we 
prove the local existence and uniqueness of geodesics for a strongly 
pseudoconvex weakly K\"ahler metric, without assuming the strong convexity of 
the metric.
 
Finally, in section~8 we deal with strongly pseudoconvex Finsler metrics of 
constant holomorphic curvature, providing a first step toward their 
classification. As a consequence of results of this section and of~[AP2] 
we get for example the following:
 
\newthm Theorem \zuno: Let $F\colon T^{1,0}M\to\R^+$ be a complete 
strongly pseudoconvex Finsler metric on a simply connected complex
manifold~$M$. Assume that 
{\smallskip
\itm{(i)} $F$ is K\"ahler;
\itm{(ii)} $F$ has constant holomorphic curvature~$-4$;
\itm{(iii)} $R(H,\bar K,\chi,\smb\chi)=R(\chi,\bar K,H,\smb\chi)$ for all 
$H$,~$K\in\ca H$, where $R$~is the curvature operator of the Chern connection;
\itm{(iv)} the indicatrices $I_F(p)=\{v\in T^{1,0}_pM\mid F(v)<1\}$ of $F$ are 
strongly convex for all~$p\in M$.
\smallskip
\noindent Then the exponential map $\exp_p\colon T^{1,0}_pM\to M$ is a 
homeomorphism, and a smooth diffeomorphism outside the origin, for any~$p\in 
M$. Furthermore, a suitable reparametrization of~$\exp_p$ induces a foliation 
of~$M$ by isometric totally geodesic holomorphic embeddings of the unit disk 
$\Delta\subset\C$ endowed with the Poincar\'e metric. In particular, $F$~is 
the Kobayashi metric of~$M$.}
 
A version of this result also holds when the holomorphic curvature is 
identically zero; the precise statement can be found in Theorem~8.10
 
\smallsect 1. Definitions and preliminaries
 
Let $M$ be a complex manifold of complex dimension~$n$. We shall
denote by $T^{1,0}M$ the holomorphic tangent bundle of~$M$, and by $\tildeM$
the complement in $T^{1,0}M$ of the zero section. The real tangent bundle
of~$M$ will be denoted by~$T_{\R}M$, and we set as usual
$T_{\C}M=T_{\R}M\otimes\C$. 
 
A {\sl complex Finsler metric} $F$ on $M$ is an 
upper semicontinuous function $F\colon T^{1,0}M\to\R^+$ satisfying
\smallskip
\item{(i)}$G=F^2$ is smooth on $\tildeM$;
\item{(ii)}$F(p;v)>0$ for all $p\in M$ and $v\in\tildeM_p$;
\item{(iii)}$F\bigl(p;\zeta v\bigr)=|\zeta|F(p;v)$ for all $p\in M$, $v\in 
T^{1,0}_p M$ and $\zeta\in\C$.
\smallskip
\noindent We shall sistematically denote by $G$ the function $G=F^2$. Note 
that it is important to ask for the smoothness of~$G$ only on~$\tildeM$: in 
fact, it is easy to see that 
$G$ is smooth on the whole of~$T^{1,0}M$ iff 
$F$~is the norm associated to a hermitian metric. In this case, we shall say 
that $F$ {\sl comes from} a hermitian metric.
 
To start, we need a few notations and general formulas. In local coordinates,
a vector $v\in T^{1,0}_pM$ is written as
$$v=v^\alpha\left.{\de\over\de z^\alpha}\right|_p,$$
where we adopt the Einstein convention. In particular, the function~$G$ is 
locally expressed in terms of the coordinates 
$\{z^1,\ldots,z^n,v^1,\ldots,v^n\}$.
We shall denote by indices like $\alpha$, $\smb\beta$ and so on the 
derivatives with respect to the $v$-coordinates; for instance,
$$G_{\alpha\smb\beta}={\de^2 G\over\de v^\alpha\de\bar{v^\beta}}.$$
On the other hand, the derivatives with respect to the $z$-coordinates will be 
denoted by indices after a semicolon; for instance,
$$G_{;\mu\nu}={\de^2 G\over\de z^\mu\de z^\nu}\qquad\hbox{or}\qquad G_{\alpha;
	\smb\nu}={\de^2 G\over\de\bar{z^\nu}\de v^\alpha}.$$
 
For our aims, we ought to focus on a 
smaller class of Finsler metrics. A complex Finsler metric $F$ will be said 
{\sl strongly pseudoconvex} if
\smallskip
\item{(iv)}the Levi matrix $(G_{\alpha\smb\beta})$ is positive definite on 
$\tildeM$.
\smallskip
\noindent This is equivalent to requiring that all the {\sl $F$-indicatrices}
$$I_F(p)=\{v\in T^{1,0}_pM\mid F(v)<1\}$$
are strongly pseudoconvexes. As we shall see in section~2, this hypothesis will 
allow us to define a hermitian metric on a suitable vector bundle. 
 
The main (actually, almost the unique) property of the function $G$ is its 
(1,1)-hom\-ogen\-eity: we have
$$G(p;\zeta v)=\zeta\smb\zeta\,G(p;v)\neweq\eqdhom$$
for all $(p;v)\in T^{1,0}M$ and $\zeta\in\C$. We now collect a number of 
formulas we shall use later on which are consequences of \eqdhom. First of 
all, differentiating with respect to $v^\alpha$ and $\bar{v^\beta}$ we get
$$\eqalign{G_\alpha(p;\zeta v)&=\smb\zeta G_\alpha(p;v),\cr
	G_{\alpha\smb\beta}(p;\zeta v)&=G_{\alpha\smb\beta}(p;v),\cr
	G_{\alpha\beta}(p;\zeta v)&=(\smb\zeta/\zeta)G_{\alpha\beta}
	(p;v).\cr}\neweq\eqdunomez$$
Thus differentiating with respect to $\zeta$ or $\smb\zeta$ and then 
setting $\zeta=1$ we get
$$G_{\alpha\smb\beta}\,\bar{v^\beta}=G_\alpha,\qquad
	G_{\alpha\beta}\,v^\beta=0,\neweq\eqdf$$
and
$$G_{\alpha\beta\gamma}\,v^\gamma=-G_{\alpha\beta},\qquad
	G_{\alpha\beta\smb\gamma}\,\bar{v^\gamma}=G_{\alpha\beta},\qquad
	G_{\alpha\smb\beta\gamma}\,v^\gamma=0,\neweq\eqde$$
where everything is evaluated at $(p;v)$.
 
On the other hand, differentiating directly \eqdhom\ with respect to $\zeta$ 
or $\smb\zeta$ and putting eventually $\zeta=1$ we get
$$G_\alpha\,v^\alpha=G,\qquad G_{\alpha\beta}\,v^\alpha v^\beta=0,\qquad
	G_{\alpha\smb\beta}\,v^\alpha\bar{v^\beta}=G.\neweq\eqdb$$
It is clear that we may get other formulas applying any differential operator 
acting only on the $z$-coordinates, or just by conjugation. For instance, we 
get
$$G_{\smb\alpha;\mu}\,\bar{v^\alpha}=G_{;\mu},\neweq\eqdo$$
and so on.
 
Assuming from now on (unless explicitely noted otherwise) $F$ strongly 
pseudoconvex, we get another bunch of formulas. As usual in hermitian geometry,
we shall denote by $(G^{\smb\beta\alpha})$ the 
inverse matrix of $(G_{\alpha\smb\beta})$, and we shall use it to raise indices.
 
First of all, applying $G^{\smb\beta\alpha}$ to the first
equation in \eqdf\ we get 
$$G^{\smb\beta\alpha}G_\alpha=\bar{v^\beta},\neweq\eqda$$
and thus, applying \eqdo,
$$G_{\smb\beta;\mu}G^{\smb\beta\alpha}G_{\alpha}=G_{;\mu}.\neweq\eqdd$$
Recalling that $(G^{\smb\beta\alpha})$ is the inverse matrix of $(G_{\alpha
\smb\beta})$, we may also compute derivatives of~$G^{\smb\beta\alpha}$:
$$DG^{\smb\beta\alpha}=-G^{\smb\nu\alpha}G^{\smb\beta\mu}(DG_{\mu\smb\nu}),
	\neweq\eqdc$$
where $D$ denotes any first order linear differential operator. As a 
consequence of \eqde\ and~\eqdc\ we get
$$G^{\smb\beta\alpha}_{\smb\sigma}\,\bar{v^\sigma}=-G^{\smb\nu\alpha}G^{\smb
	\beta\mu}G_{\mu\smb\nu\smb\sigma}\,\bar{v^\sigma}=0,\neweq\eqdg$$
and recalling also \eqda\ we obtain
$$G_{\smb\beta}G^{\smb\beta\alpha}_\gamma=-G_{\smb\beta}G^{\smb\beta\mu}
	G^{\smb\nu\alpha}G_{\mu\smb\nu\gamma}=-G^{\smb\nu\alpha}G_{\mu\smb\nu
	\gamma}v^\mu=0.\neweq\eqdh$$
 
\smallsect 2. The Chern-Finsler connection
 
To any hermitian metric is associated a
unique (1,0)-connection such that the metric tensor is parallel: the 
Chern connection. The main goal of this section is to define the analogue for
strongly pseudoconvex Finsler metrics. 
 
Let $\pi\colon\tildeM\to M$ denote the restriction of the canonical projection 
of $T^{1,0}M$ onto~$M$. The {\sl vertical bundle} $\ca V\subset T^{1,0}M$ is, 
by definition, the kernel of the differential $d\pi\colon T^{1,0}\tildeM\to
T^{1,0}M$. It is easy to check that $\ca V$ is a complex vector bundle of rank 
$n$ over~$\tildeM$; a local frame for $\ca V$ is given by $\{\dot\de_1,\ldots,
\dot\de_n\}$, where we set
$$\dot\de_\alpha={\de\over\de v^\alpha}\qquad\hbox{and}\qquad\de_\mu={\de\over
	\de z^\mu},$$
for $\alpha$,~$\mu=1,\ldots,n$. We shall denote by $\ca X(\ca V)$ the space of 
smooth sections of $\ca V$; more generally, $\ca X(E)$ will denote the space 
of smooth sections of any vector bundle~$p\colon E\to B$.
 
Let $j_p\colon T^{1,0}_pM\hookrightarrow T^{1,0}M$ be the inclusion and, for 
$v\in\tildeM_p$, let $k_v\colon T^{1,0}_p M\to T^{1,0}_v(T^{1,0}_pM)$ denote 
the usual identification. Then we get a natural isomorphism
$$\iota_v=d(j_{\pi(v)})_v\circ k_v\colon T^{1,0}_{\pi(v)}M\to\ca V_v,$$
and, by restriction, the all-important natural section $\iota\colon\tildeM\to
\ca V$ given by
$$\iota(v)=\iota_v(v)\in\ca V_v.$$
In local coordinates,
$$\iota_v\left(\left.{\de\over\de z^\alpha}\right|_{\pi(v)}\right)=
	\dot\de_\alpha|_v;$$
in particular, if $v=v^\alpha(\de/\de z^\alpha)$ then
$$\iota(v)=v^\alpha\dot\de_\alpha|_v.$$
$\iota$ is called the {\sl radial vertical vector field.}
 
The first observation is that a strongly pseudoconvex Finsler metric $F$ 
defines a hermitian metric on the vertical 
bundle~$\ca V$. Indeed, if $v\in\tildeM$ and $W_1$, $W_2\in\ca V_v$, with 
$W_j=W^\alpha_j\dot\de_\alpha$, we set
$$\langle W_1,W_2\rangle_v=G_{\alpha\smb\beta}(v)W^\alpha_1\bar{W^\beta_2}.$$
Being $F$ strongly pseudoconvex, $\langle\,,\rangle$ is a hermitian metric. 
Note that the third equation in \eqdb\ says that
$$G=\langle\iota,\iota\rangle;$$
so $\iota$ is an isometric embedding of~$\tildeM$ into~$\ca V$.
 
Following Kobayashi [K3], we now consider the Chern connection~$D$ on the 
vector bundle~$\ca V$: it is the unique (1,0)-connection on~$\ca V$ such that 
the hermitian structure previously defined is parallel. In other words, 
$D\colon\ca X(\ca V)\to\ca X(T^*_{\C}\tildeM\otimes\ca V)$ is such that
$$X\langle V,W\rangle=\langle\nablab_XV,W\rangle+\langle V,\nablab_{\bar X}W
	\rangle,$$
for any $X\in T^{1,0}\tildeM$ and $V$,~$W\in\ca X(\ca V)$.
 
In local coordinates, the connection matrix $(\omega^\alpha_\beta)$ is given 
by
$$\omega^\alpha_\beta=G^{\smb\tau\alpha}\de G_{\beta\smb\tau}=
	\tilde\Gamma^\alpha_{\beta;\mu}\,dz^\mu+\tilde\Gamma^\alpha
	_{\beta\gamma}\,dv^\gamma,$$
where
$$\tilde\Gamma^\alpha_{\beta\gamma}=G^{\smb\tau\alpha}G_{\beta\smb\tau\gamma}
	\qquad\hbox{and}\qquad\tilde\Gamma^\alpha_{\beta;\mu}=G^{\smb\tau\alpha}
	G_{\beta\smb\tau;\mu}.$$
 
This is only part of the connection we are looking for: our next goal is to 
canonically extend~$D$ to a (1,0)-connection on $T^{1,0}\tildeM$. Let us 
consider the bundle map $\Lambda\colon T^{1,0}\tildeM\to\ca V$ defined by
$$\Lambda(X)=\nablab_X\iota,$$
and set $\ca H=\ker\Lambda\subset T^{1,0}\tildeM$. We claim that $\ca H$ is a
{\sl horizontal bundle,} that is $T^{1,0}\tildeM=\ca H\oplus\ca V$. Indeed, in 
local coordinates
$$\Lambda(X)=[\dot X^\alpha+\omega^\alpha_\beta(X)v^\beta]\dot\de_\alpha,$$
where $X=X^\mu\de_\mu+\dot X^\alpha\dot\de_\alpha$. Then a local frame for
$\ca H$ is given by $\{\delta_1,\ldots,\delta_n\}$, where
$$\delta_\mu=\de_\mu-\tilde\Gamma^\alpha_{\beta;\mu}v^\beta\dot\de_\alpha$$
--- note that $\tilde\Gamma^\alpha_{\beta\gamma}v^\beta\equiv0$ --- and the 
claim is proved.
 
It is not difficult to check (see [AP4] for a coordinate-free proof) that 
setting
$$\Theta(\dot\de_\alpha)=\delta_\alpha$$
for $\alpha=1,\ldots,n$ we get a well-defined global bundle isomorphism
$\Theta\colon\ca V\to\ca H$; then we can define a (1,0)-connnection~$D$ 
on~$\ca H$ just by setting
$$\nablab_X H=\Theta\bigl[\nablab_X(\Theta^{-1}H)\bigr]$$
for any $X\in T_\C\tildeM$ and $H\in\ca X(\ca H)$. By linearity, this yields a 
(1,0)-connection on $T^{1,0}\tildeM$, still denoted by~$D$: the {\sl 
Chern-Finsler connection.}
 
Using the bundle isomorphism $\Theta\colon\ca V\to\ca H$ we can also transfer
the hermitian structure~$\langle\,,\rangle$ on~$\ca H$ just by setting 
$$\forevery{H,K\in\ca H_v}\langle H,K\rangle_v=\langle\Theta^{-1}(H),
	\Theta^{-1}(K)\rangle_v,$$
and then we can define a hermitian structure on $T^{1,0}\tildeM$ by requiring 
$\ca H$ be orthogonal to~$\ca V$. It is easy to check then that $D$ is the 
Chern connection associated to this hermitian structure, that is 
$$X\langle Y,Z\rangle=\langle\nablab_XY,Z\rangle+\langle Y,\nabla_{\bar X}Z
	\rangle$$
for any $X\in T^{1,0}\tildeM$ and $Y$,~$Z\in\ca X(T^{1,0}\tildeM)$.
 
>From now on we shall work only with the frame
$\{\delta_\mu,\dot\de_\alpha\}$ and its dual co-frame $\{dz^\mu,\psi^\alpha\}$
given by
$$\psi^\alpha=dv^\alpha+\Gamma^\alpha_{;\mu}\,dz^\mu=dv^\alpha+G^{\smb\tau
	\alpha}G_{\smb\tau;\mu}\,dz^\mu,$$
where we have set
$$\Gamma^\alpha_{;\mu}=\tilde\Gamma^\alpha_{\beta;\mu}v^\beta=G^{\smb\tau
	\alpha}G_{\smb\tau;\mu}.$$
Writing
$$\omega^\alpha_\beta=\Gamma^\alpha_{\beta;\mu}\,dz^\mu+\Gamma^\alpha_{\beta
	\gamma}\psi^\gamma,$$
we get
$$\eqalign{\Gamma^\alpha_{\beta\gamma}&=G^{\smb\tau\alpha}G_{\beta\smb\tau
	\gamma}=\Gamma^\alpha_{\gamma\beta},\cr
	\Gamma^\alpha_{\beta;\mu}&=G^{\smb\tau\alpha}\delta_\mu(G_{\beta\smb
	\tau})=G^{\smb\tau\alpha}(G_{\beta\smb\tau;\mu}
	-G_{\beta\smb\tau\gamma}\Gamma^\gamma_\mu).\cr}\neweq\eqdgamma$$
Note that
$$\Gamma^\alpha_{\beta;\mu}=\dot\de_\beta(\Gamma^\alpha_{;\mu})\qquad
	\hbox{and}\qquad\Gamma^\alpha_{;\mu}=\Gamma^\alpha_{\beta;\mu}
	v^\beta;\neweq\eqdvt$$
in particular, this is exactly the connection introduced by Rund~[Ru2].
 
So we have described a canonical splitting of the holomorphic 
tangent bundle of~$\tildeM$ in a vertical and a horizontal bundle, and defined 
a canonical connection on it, preserving this splitting. In the following 
subsections we shall begin the study of this connection, introducing torsions 
and curvatures; here we first describe a few properties of the splitting.
 
First of all, the next lemma shows that the local frames $\{\delta_1,\ldots,
\delta_n\}$ enjoy some nice and convenient properties:
 
\newthm Lemma \dtre: Let $D$ be the Chern-Finsler connection 
associated to a strongly pseudoconvex Finsler metric $F$, and let $\{\delta_1,
\ldots,\delta_n\}$ be the corresponding local horizontal frame. Then
{\smallskip
\itm{(i)} $[\delta_\mu,\delta_\nu]=0$ for all $1\le \mu,\nu\le n$;
\itm{(ii)} $[\delta_\mu,\dot\de_\alpha]=\Gamma^\sigma_{\alpha;\mu}\dot
\de_\sigma$ for 
all $1\le\alpha,\mu\le n$;
\itm{(iii)} $\delta_\mu(G)=\delta_{\smb\mu}(G)=0$ for all $1\le\mu\le n$;
\itm{(iv)} $\delta_{\smb\mu}(G_\alpha)=0$ for all $1\le\alpha,\mu\le n$.}
 
\pf (i) If suffices to compute. First of all,
$$[\delta_\mu,\delta_\nu]=(\Gamma^\alpha_{;\mu\nu}-\Gamma^\alpha_{;\nu\mu}+
	\Gamma^\alpha_{\sigma;\nu}\Gamma^\sigma_{;\mu}-\Gamma^\alpha_{\sigma;
	\mu}\Gamma^\sigma_{;\nu})\dot\de_\alpha,$$
where $\Gamma^\alpha_{;\mu\nu}=\de_\nu(\Gamma^\alpha_{;\mu})$ and so on. Now,
$$\eqalign{\Gamma^\alpha_{;\mu\nu}&=G^{\smb\tau\alpha}(G_{\smb\tau;\mu\nu}-
	G_{\sigma\smb\tau;\nu}\Gamma^\sigma_{;\mu}),\cr
	\Gamma^\alpha_{;\nu\mu}&=G^{\smb\tau\alpha}(G_{\smb\tau;\nu\mu}-
	G_{\sigma\smb\tau;\mu}\Gamma^\sigma_{;\nu}),\cr}\qquad
\eqalign{\Gamma^\alpha_{\sigma;\nu}\Gamma^\sigma_{;\mu}&=G^{\smb\tau\alpha}
	(G_{\sigma\smb\tau;\nu}\Gamma^\sigma_{;\mu}-
	G_{\sigma\smb\tau\rho}\Gamma^\rho_{;\nu}\Gamma^\sigma_{;\mu}),\cr
	\Gamma^\alpha_{\sigma;\mu}\Gamma^\sigma_{;\nu}&=G^{\smb\tau\alpha}
	(G_{\sigma\smb\tau;\mu}\Gamma^\sigma_{;\nu}-
	G_{\sigma\smb\tau\rho}\Gamma^\rho_{;\mu}\Gamma^\sigma_{;\nu}),\cr}$$
and the assertion follows. Note that we have actually proved that
$$\delta_\nu(\Gamma^\alpha_{;\mu})=\delta_\mu(\Gamma^\alpha_{;\nu}).\neweq
	\eqdcom$$
 
(ii) Indeed,
$$[\delta_\mu,\dot\de_\alpha]=[\de_\mu-\Gamma^\sigma_{;\mu}\dot\de_\sigma,
	\dot\de_\alpha]=\dot\de_\alpha(\Gamma^\sigma_{;\mu})\dot\de_\sigma=
	\Gamma^\sigma_{\alpha;\mu}\dot\de_\sigma.$$
 
(iii) In fact, using \eqdd\ we get
$$\delta_\mu(G)=G_{;\mu}-\Gamma^\sigma_{;\mu}G_\sigma=G_{;\mu}-
	G^{\smb\tau\sigma}G_{\smb\tau;\mu}G_\sigma=G_{;\mu}-G_{;\mu}=0.$$
 
(iv) Finally,
$$\delta_{\smb\mu}(G_\alpha)=G_{\alpha;\smb\mu}-\Gamma^{\smb\tau}_{;\smb\mu}
	G_{\alpha\smb\tau}=G_{\alpha;\smb\mu}-G_{\alpha;\smb\mu}=0,$$
where $\Gamma^{\smb\tau}_{;\smb\mu}=\bar{\Gamma^\tau_{;\mu}}$.\qedn
 
The philosophical idea behind our work is that to study the geometry of a
complex Finsler metric one should transfer everything (or most of it) in the
horizontal bundle, and then apply the usual techniques of hermitian geometry
there. We shall better substantiate this idea later, for instance in
sections~6 and~7 discussing variation formulas; here we begin to show how to
lift objects (e.g., vector fields) from the tangent bundle up to~$\ca H$. 
 
The main tool is provided by the horizontal analogues of the 
isomorphisms~$\iota_v$. If $v\in\tildeM$, we set
$$\chi_v=\Theta_v\circ\iota_v\colon T_{\pi(v)}M\to\ca H_v.$$
The {\sl horizontal radial vector field} $\chi\in\ca X(\ca H)$ is then defined 
by
$$\chi=\Theta\circ\iota;$$
in local coordinates, if $v=v^\alpha(\de/\de z^\alpha)|_p$ we have
$$\chi(v)=v^\alpha\delta_\alpha|_v.$$
 
Using the isomorphisms~$\chi_v$ we can induce an embedding of~$\tildeM$ 
into~$\ca H$ which respects the Lie algebra structure. To be precise,
a vector field $\xi\in\ca X(T^{1,0}M)$ may be lifted in two 
different ways to vector fields in $T^{1,0}\tildeM$: via the {\sl horizontal 
lift}
$$\xi^H(v)=\chi_v\Bigl(\xi\bigl(\pi(v)\bigr)\Bigr),$$
and via the {\sl vertical lift}
$$\xi^V(v)=\iota_v\Bigl(\xi\bigl(\pi(v)\bigr)\Bigr).$$
A consequence of Lemma \dtre\ is that the horizontal lift is a Lie algebra 
homomorphism:
 
\newthm Proposition \dquattro: Let $D$ be the Chern-Finsler connection 
associated to a strongly pseudoconvex Finsler metric~$F$ on a complex 
manifold~$M$. Then:
{\smallskip
\itm{(i)}$[\ca X(\ca H),\ca X(\ca H)]\subset\ca X(\ca H)$ and $[\ca X(\ca 
V),\ca X(\ca V)]\subset\ca X(\ca V)$;
\itm{(ii)}if $\xi_1$,~$\xi_2\in\ca X(\tildeM)$ then $[\xi^H_1,\xi^H_2]=[\xi_1,
\xi_2]^H$, $[\xi^V_1,\xi^V_2]=0$ and $[\xi_1^H,\xi_2^V]\in\ca X(\ca V)$.}
 
\pf (i) Take $H_1$, $H_2\in\ca X(\ca H)$. Locally, $H_j=H^\mu_j\delta_\mu$; 
hence
$$[H_1,H_2]=\bigl(H^\nu_1\delta_\nu(H^\mu_2)-H^\nu_2\delta_\nu(H^\mu_1)\bigr)
	\delta_\mu
	\neweq\eqdunoz$$
(where we used Lemma \dtre) is horizontal. Analogously, if $V_1$, $V_
2\in\ca X(\ca V)$ with $V_j=V^\alpha_j\dot\de_\alpha$, we get
$$[V_1,V_2]=\bigl(V^\beta_1\dot\de_\beta(V^\alpha_2)-V^\beta_2\dot\de_\beta
	(V^\alpha_1)\bigr)\dot\de_\alpha,\neweq\eqdduez$$
which is vertical.
 
(ii) Locally, $\xi_j=\xi^\mu_j(\de/\de z^\mu)$ and $\xi^H_j=(\xi^\mu_j\circ\pi)
\delta_\mu$; so \eqdunoz\ yields
$$[\xi^H_1,\xi^H_2]=\bigl((\xi^\nu_1\circ\pi)\delta_\nu(\xi^\mu_2\circ\pi)-
	(\xi^\nu_2\circ\pi)\delta_\nu(\xi^\mu_1\circ\pi)\bigr)\delta_\mu.$$
Now $\delta_\nu(\xi^\mu_j\circ\pi)=(\de\xi^\mu_j/\de z^\nu)\circ\pi$; therefore
$$[\xi^H_1,\xi^H_2]=\left[\left(\xi^\nu_1{\de\xi^\mu_2\over\de z^\nu}-\xi^\nu_2
	{\de\xi^\mu_1\over\de z^\nu}\right)\circ\pi\right]\delta_\mu=
	[\xi_1,\xi_2]^H.$$
On the other hand, $\xi^V_j=(\xi^\alpha_j\circ\pi)\dot\de_\alpha$ and
$\dot\de_\beta(\xi^\alpha_j\circ\pi)=0$ yield
$$[\xi^V_1,\xi^V_2]=0.$$
Finally,
$$[\xi^H_1,\xi^V_2]=\left[\left(\xi^\mu_1{\de\xi^\alpha_2\over\de z^\mu}\right)
	\circ\pi+\bigl((\xi^\mu_1\xi^\beta_2)\circ\pi\bigr)\Gamma^\alpha_
	{\beta;\mu}\right]\dot\de_\alpha,$$
again by Lemma~\dtre.\qedn
 
Note that, as a consequence of (ii), the obvious map of $\ca X(\ca 
V)$ into $\ca X(\ca H)$ induced by the complex horizontal map $\Theta\colon\ca 
V\to\ca X$ is {\it not} an isomorphism of Lie algebras; it suffices 
to remark that $\Theta(\xi^V)=\xi^H$ for all $\xi\in\ca X(\tildeM)$.
 
\smallsect 3. Torsions and k\"ahlerianity
 
As it may be expected, the next step is the study of the Chern-Finsler 
connection is to describe its torsion(s) and clarify their geometrical 
meaning.
 
The tangent bundle $T^{1,0}M$ (and hence $\tildeM$ too) is naturally equipped 
with a $T^{1,0}\tildeM$-valued global (1,0)-form, the {\sl canonical form}
$$\eta=dz^\mu\otimes\de_\mu+dv^\alpha\otimes\dot\de_\alpha\in\ca X(\Wedge^{1,0}
	\tildeM\otimes T^{1,0}\tildeM).$$
It is easy to see that as soon as we have a strongly pseudoconvex 
Finsler metric --- and hence the canonical splitting $T^{1,0}\tildeM=\ca H
\oplus\ca V$ --- one has
$$\eta=dz^\mu\otimes\delta_\mu+\psi^\alpha\otimes\dot\de_\alpha.$$
 
Extending as usual the Chern-Finsler connection~$D$ to an exterior 
differential (still denoted by~$D$) on $T^{1,0}\tildeM$-valued 
differential forms, it is very natural to consider the {\sl torsion}~$D\eta$ 
of the connection. Since $\eta$~is a (1,0)-form, $D\eta$ splits in the sum of 
a (2,0)-form~$\theta$ and a (1,1)-form~$\tau$. We shall call $\theta$ the {\sl
$(2,0)$-torsion} of the Chern-Finsler connection, and $\tau$ the {\sl 
$(1,1)$-torsion} of the Chern-Finsler connection.
 
Locally, we may write
$$\theta=\theta^\mu\otimes\delta_\mu+\dot\theta^\alpha\otimes\dot\de_\alpha
	\qquad\hbox{and}\qquad\tau=\tau^\alpha\otimes\dot\de_\alpha,$$
where, setting $\Gamma^\alpha_{\smb\beta;\mu}=\dot\de_{\smb\beta}(\Gamma^\alpha
_{;\mu})$,
$$\eqalign{\tau^\alpha&=\smb\de\psi^\alpha=-\delta_{\smb\nu}(\Gamma^\alpha
	_{;\mu})\,dz^\mu\wedge d\smb z^\nu-\Gamma^\alpha_{\smb\beta;\mu}
	\,dz^\mu\wedge\bar{\psi^\beta};\cr
	\theta^\mu&=-dz^\nu\wedge\omega^\mu_\nu={\textstyle{1\over2}}
	[\Gamma^\mu_{\nu;\sigma}-\Gamma^\mu_{\sigma;\nu}]\,dz^\sigma
	\wedge dz^\nu+\Gamma^\mu_{\nu\gamma}\,\psi^\gamma\wedge dz^\nu;\cr}
	\neweq\eqTau$$
and
$$\eqalign{\dot\theta^\alpha&=\de\psi^\alpha-\psi^\beta\wedge\omega^\alpha_\beta
	\cr
	&={\textstyle{1\over2}}[\delta_\mu(\Gamma^\alpha_{;\nu})-\delta_\nu
	(\Gamma^\alpha_{;\mu})]\,dz^\mu\wedge dz^\nu+[\dot\de_\beta(\Gamma
	^\alpha_{;\mu})-\Gamma^\alpha_{\beta;\mu}]\,\psi^\beta\wedge dz^\mu
	+{\textstyle{1\over2}}[\Gamma^\alpha_{\beta\gamma}-\Gamma^\alpha
	_{\gamma\beta}]\,\psi^\beta\wedge\psi^\gamma\cr
	&=0,\cr}\neweq\eqthp$$
by \eqdcom, \eqdgamma\ and \eqdvt.
 
One may wonder whether these torsions are the right generalizations of the 
usual torsion in the hermitian case. The answer is a double yes. First of all, 
a standard argument using the definitions shows that torsions and covariant 
derivative are related as usual:
$$\eqalign{\nablab_XY-\nablab_YX&=[X,Y]+\theta(X,Y),\cr
	\nablab_X\bar{Y}-\nabla_{\bar Y}X&=[X,\bar{Y}]+\tau(X,\bar Y)
	+\bar\tau(X,\bar Y),\cr}\neweq\eqtor$$
for any $X$,~$Y\in\ca X(T^{1,0}\tildeM)$, where, by definition, 
$$\nablab_X\bar Y=\bar{\nabla_{\bar X}Y}.$$
 
Furthermore, the vanishing of (part of) the (2,0)-torsion can be again 
interpreted as a K\"ahler condition --- but with some care, because $\theta$ 
is composed by a horizontal part and a mixed part.  
To be precise, we shall say that a differential form~$\gamma$ on~$\tildeM$ is 
{\sl horizontal} if
it vanishes contracted with any $V\in\ca X(\ca V)$. The 
decomposition $T^{1,0}\tildeM=\ca H\oplus\ca V$ induces a projection $p^*_H$ of 
the differential forms onto the horizontal forms; the {\sl horizontal part} 
of a form~$\gamma$ is then~$p^*_H(\gamma)$.
 
There is a corresponding projection on the vertical forms, of course, but we
shall not need it now because the vertical part of both torsions $\theta$ and 
$\tau$ is zero. For this reason, the form $\theta-p^*_H(\theta)$ will be called
the {\sl mixed part} of $\theta$. In local coordinates,
$$p^*_H(\theta)=(\Gamma^\sigma_{\nu;\mu}\,dz^\mu\wedge dz^\nu)\otimes\delta_
	\sigma\qquad\hbox{and}\qquad \theta-p^*_H(\theta)=(\Gamma^\sigma_{\nu
	\gamma}\,\psi^\gamma\wedge dz^\nu)\otimes\delta_\sigma.$$
The next proposition discusses the meaning of the vanishing of the 
(2,0)-torsion~$\theta$ or of one of its parts.
 
\newthm Proposition \tvan: Let $F$ be a strongly pseudoconvex Finsler metric 
on a complex manifold~$M$. Then:
{\smallskip
\itm{(i)} the mixed part of the $(2,0)$-torsion vanishes iff $F$ comes from a 
hermitian metric;
\itm{(ii)} $\theta$ vanishes iff $F$ comes from a hermitian K\"ahler metric.}
 
\pf (i) The mixed part of the torsion vanishes iff $G_{\beta\smbar\mu\gamma}=
0$ for all $\beta$, $\mu$ and $\gamma$. Conjugating, this is equivalent to 
having $\dot\de_\gamma(G_{\beta\smbar\mu})=\dot\de_{\smbar\gamma}(G_{\beta
\smbar\mu})=0$, that is $G_{\beta\smbar\mu}(v)$ depends only on $\pi(v)$ --- 
and this happens iff $F$ comes from a hermitian metric.
 
(ii) It follows from (i) and the fact that when $F$ comes from a hermitian 
metric $g=(g_{\alpha\smbar\beta})$ one has
$$\Gamma^\alpha_{\beta;\mu}=g^{\smbar\tau\alpha}{\de g_{\beta\smbar\tau}\over
	\de z^\mu}.$$
\qedn
 
For this reason we say that a strongly pseudoconvex Finsler metric $F$ is {\sl 
strongly K\"ahler} if the horizontal part of the (2,0)-torsion vanishes, that 
is iff
$$\forevery{H,K\in\ca H}\theta(H,K)=0.$$
This is exactly the notion of k\"ahlerianity introduced by Rund~[Ru2]. 
However, as we shall see later on (see sections~6 and~7), studying the 
geometry of a strongly pseudoconvex Finsler metric it turns out that this 
assumption is too strong and not quite natural. So it is appropriate to 
introduce two more notions of k\"ahlerianity. We shall say that $F$~is 
{\sl K\"ahler} if
$$\forevery{H\in\ca H}\theta(H,\chi)=0,$$
and that $F$ is {\sl weakly K\"ahler} if
$$\forevery{H\in\ca H}\langle \theta(H,\chi),\chi\rangle=0.$$
In local coordinates, $F$ is strongly K\"ahler iff
$$\Gamma^\alpha_{\mu;\nu}=\Gamma^\alpha_{\nu;\mu};$$
it is K\"ahler iff
$$\Gamma^\alpha_{\mu;\nu}v^\mu=\Gamma^\alpha_{\nu;\mu}v^\mu;$$
it is weakly K\"ahler iff
$$G_\alpha[\Gamma^\alpha_{\mu;\nu}-\Gamma^\alpha_{\nu;\mu}]v^\mu=0,$$
that is iff
$$0=[G_{\mu;\nu}-G_{\nu;\mu}+G_{\nu\sigma}\Gamma^\sigma_{;\mu}]v^\mu
	=[G_{\mu\smb\tau;\nu}-G_{\nu\smb\tau;\mu}+G_{\nu\sigma\smb\tau}
	\Gamma^\sigma_{;\mu}]v^\mu\bar{v^\tau}.$$
In particular, if $F$ comes from a hermitian metric then these three 
conditions are all equivalent to the usual K\"ahler condition, because $G_{\nu
\sigma\smb\tau}\equiv0$ for a Finsler metric coming from a hermitian metric. 
 
There are other characterizations of strongly K\"ahler Finsler metrics. To $F$ 
we may associate the {\sl fundamental form}
$$\Phi=iG_{\alpha\smbar\beta}\,dz^\alpha\wedge d\bar{z^\beta},$$
which is a well-defined real (1,1)-form on $\tildeM$.
Then the strong K\"ahler condition is equivalent to the vanishing of 
the horizontal part of $d\Phi$. To express it more clearly, set 
$$d_H=p^*_H\circ d,\qquad \de_H=p^*_H\circ\de\qquad\hbox{and}\qquad\smb
	\de_H=p^*_H\circ\smbar\de,$$ 
so that again $d_H=\de_H+\smbar\de_H$. 
 
\newthm Theorem \tkahler: Let $F$ be a strongly pseudoconvex Finsler metric 
on a complex manifold~$M$. Then the following assertions are equivalent:
{\smallskip
\itm{(i)} $F$ is a strongly K\"ahler Finsler metric;
\itm{(ii)} $\nablab_HK-\nablab_KH=[H,K]$ for all $H$,~$K\in\ca X(\ca H)$;
\itm{(iii)}$d_H\Phi=0$;
\itm{(iv)}$\de_H\Phi=0$;
\itm{(v)} for any $p_0\in M$ there is a neighbourhood $U$ of $p_0$ in $M$ and 
a real-valued function $\phi\in C^\infty\bigl(\pi^{-1}(U)\bigr)$ such that 
$\Phi=i\de_H\smbar\de_H\phi$ on $\pi^{-1}(U)$.}
 
\pf (i)${}\Longleftrightarrow{}$(ii) follows from \eqtor. 
 
(iii)${}\Longleftrightarrow{}$(iv) holds simply because $\Phi$ is a real 
(1,1)-form.
 
(iv)${}\Longleftrightarrow{}$(i). Indeed, \eqdgamma\ yields
$$\de\Phi(X,Y,\bar Z)=i\langle \theta(X,Y),Z\rangle$$
for all $X$,~$Y$,~$Z\in T^{1,0}\tildeM$;
hence $\de_H\Phi$ vanishes iff $p^*_H\circ \theta$ vanishes, that is iff $F$ is 
strongly K\"ahler.
 
(v)${}\Longrightarrow{}$(iv) follows from Lemma~\dtre.(i).
 
(iii)${}\Longrightarrow{}$(v). 
Let $\gamma$ be any horizontal form. In local coordinates, defined on a 
coordinate neighbourhood of the form $\pi^{-1}(U)$, one has
$$\gamma|_{(p;v)}=\gamma_{A\smbar B}(p;v)\,dz^A\wedge d\bar{z^B},$$
for suitable multi-indices $A$ and $\smbar B$. On $U$ we may then consider the 
family of forms
$$\gamma_v|_p=\gamma_{A\smbar B}(p;v)\,dz^A\wedge d\bar{z^B},$$
where here $\{dz^j\}$ is the dual frame of $\{\de/\de z^j\}$; in other words, 
we are considering the $v$-coordinates just as parameters. 
 
The gist is that the following formula holds:
$$(d_H\gamma)_v=d(\gamma_v).$$
Then we may now apply the Dolbeault and Serre theorems (with parameters) to 
$\Phi_v$ in a possibly smaller neighbourhood of $p_0$ --- still denoted by $U$ 
--- to get a function $\phi_v\in C^\infty(U,\R)$ depending smoothly on $v$ 
such that $\Phi_v=i\de\smbar\de\phi_v$. Then setting 
$$\phi(p;v)=\phi_v(p)$$ 
we get $\Phi=i\de_H\smbar\de_H\phi$, as required.\qedn
 
>From this point of view, a strongly pseudoconvex Finsler metric is 
K\"ahler iff
$$d_H\Phi(\cdot,\chi,\cdot)\equiv0,$$
and it is weakly K\"ahler iff
$$d_H\Phi(\cdot,\chi,\smb\chi)\equiv0.$$
 
We end this section pointing out that also the vanishing of the (1,1)-torsion 
$\tau$ has a nice geometric meaning:
 
\newthm Proposition \imm: The $(1,1)$-torsion $\tau$ vanishes iff the frame
$\{\delta_\mu,\dot\de_\alpha\}$ is holomorhic.
 
\pf Indeed the frame $\{\delta_\mu,\dot\de_\alpha\}$ is holomorphic iff its 
dual coframe $\{dz^\mu,\psi^\alpha\}$ is, which happens iff the forms 
$\psi^\alpha$ are holomorphic, that is iff $\tau^\alpha=\smb\de\psi^\alpha=0$ 
for $\alpha=1,\ldots,n$.\qedn
 
\smallsect 4. The curvature tensor
 
The {\sl curvature tensor} $R\colon\ca 
X(T^{1,0}\tildeM)\to\ca X(\Wedge^2(T^*_\C\tildeM)\otimes T^{1,0}\tildeM)$ of 
the Chern-Finsler connection is given 
by $R=D\circ D$, that is
$$\forevery{X\in\ca X(T^{1,0}\tildeM)} R_X=D(DX).$$
Analogously we have the {\sl curvature operator} $\Omega\in\ca X(\Wedge^2
(T^*_\C\tildeM)\otimes\Wedge^{1,0}\tildeM\otimes T^{1,0}\tildeM)$ defined by
(cf.~also~[K3])
$$\Omega(X,Y)Z=R_Z(X,Y).$$
Locally, $\Omega$ is given by
$$\Omega=\Omega^\alpha_\beta\otimes[dz^\beta\otimes\delta_\alpha+\psi^\beta
	\otimes\dot\de_\alpha],$$
where
$$\Omega^\alpha_\beta=d\omega^\alpha_\beta-\omega^\gamma_\beta\wedge\omega
	^\alpha_\gamma.$$
 
Decomposing $\Omega$ into types, we  get 
$$\Omega=\Omega'+\Omega'',$$
where $\Omega'$ is a (2,0)-form and $\Omega''$ a (1,1)-form. Locally,
$$(\Omega')^\alpha_\beta=\de\omega^\alpha_\beta-\omega^\gamma_\beta\wedge
	\omega^\alpha_\gamma,\qquad(\Omega'')^\alpha_\beta=\smb\de\omega^
	\alpha_\beta.$$
$\Omega$ has no (0,2)-components because the connection forms are (1,0)-forms. 
Actually, even $\Omega'$ vanishes: indeed, by definition
$$\omega^\beta_\alpha=G^{\smbar\tau\beta}\de G_{\alpha\smbar\tau}.$$
So
$$\eqalign{\de\omega^\beta_\alpha&=\de G^{\smbar\tau\beta}\wedge\de G_{\alpha
	\smbar\tau}=-G^{\smbar\tau\mu}G^{\smbar\nu\beta}\de G_{\mu\smbar\nu}
	\wedge\de G_{\alpha\smbar\tau}\cr
	&=(G^{\smbar\tau\mu}\de G_{\alpha\smbar\tau})\wedge(G^{\smbar\nu\beta}
	\de G_{\mu\smbar\nu})=\omega^\mu_\alpha\wedge\omega_\mu^\beta.\cr}
	\neweq\eqdom$$
So $\Omega=\Omega''$ and
$$\Omega^\alpha_\beta=\smb\de\omega^\alpha_\beta,$$
exactly as in the hermitian case.
 
The relation between curvature and covariant derivatives is the usual one:
$$\eqalign{\nablab_X\nablab_Y-\nablab_Y\nablab_X&=\nabla_{[X,Y]}\cr
	\nablab_X\nabla_{\bar Y}-\nabla_{\bar Y}\nablab_X&=\nabla_{[X,\bar Y]}
	+\Omega(X,\bar Y);\cr
	\nabla_{\bar X}\nabla_{\bar Y}-\nabla_{\bar Y}\nabla_{\bar X}&=
	\nabla_{[\bar X,\bar Y]},\cr}$$
for any $X$,~$Y\in\ca X(T^{1,0}\tildeM)$.
 
We can also recover the Bianchi identities in this setting:
 
\newthm Proposition \dB: Let $D\colon\ca X(T^{1,0}\tildeM)\to\ca X(T^*_\C
\tildeM\otimes T^{1,0}\tildeM)$ be the complex linear connection on~$\tildeM$ 
induced by a good complex vertical connection. Then
$$\displaylines{D\theta=\eta^H\wedge\Omega,\cr
	D\tau=\eta^V\wedge\Omega,\cr
	D\Omega=0,\cr}$$
where $\eta^H=dz^\mu\otimes\delta_\mu$ and $\eta^V=\psi^\alpha\otimes\dot\de
_\alpha$.
 
\pf It suffices to compute. First of all,
$$\eqalign{\smb\de\theta^\mu&=dz^\nu\wedge\smb\de\omega^\mu_\nu=dz^\nu\wedge
	\Omega^\mu_\nu,\cr
	\de\theta^\mu+\theta^\nu\wedge\omega^\mu_\nu&=dz^\nu\wedge\de\omega
	^\mu_\nu-dz^\nu\wedge\omega^\sigma_\nu\wedge\omega^\mu_\sigma=0,\cr}$$
by \eqdom, and so $D\theta=\eta^H\wedge\Omega$. Next
$$\displaylines{\smb\de\tau^\alpha=0,\cr
	\de\tau^\alpha+\tau^\beta\wedge\omega^\alpha_\beta=\de\smb\de
	\psi^\alpha+\smb\de\psi^\beta\wedge\omega^\alpha_\beta=-\smb\de\de
	\psi^\alpha+\smb\de\psi^\beta\wedge\omega^\alpha_\beta=\psi^\beta
	\wedge\smb\de\omega^\alpha_\beta,\cr}$$
by \eqthp, and so $D\tau=\eta^V\wedge\Omega$. Finally, $\smb\de\Omega^\alpha
_\beta=0$ and
$$\de\Omega^\alpha_\beta-\omega^\gamma_\beta\wedge\Omega^\alpha_
	\gamma+\Omega^\gamma_\beta\wedge\omega^\alpha_\gamma=\de\smb\de
	\omega^\alpha_\beta-\omega^\gamma_\beta\wedge\smb\de\omega^\alpha_
	\gamma+\smb\de\omega^\gamma_\beta\wedge\omega^\alpha_\gamma
	=-\smb\de(\de\omega^\alpha_\beta-\omega^\gamma_\alpha\wedge
	\omega^\beta_\gamma)=0,$$
by \eqdom.\qedn
 
In local coordinates, the curvature operator is given by
$$\Omega^\alpha_\beta=R^\alpha_{\beta;\mu\smbar\nu}\,dz^\mu\wedge d\smb z^\nu+
	R^\alpha_{\beta\delta;\smbar\nu}\,\psi^\delta\wedge d\smb z^\nu+
	R^\alpha_{\beta\smbar\gamma;\mu}\,dz^\mu\wedge\bar{\psi^\gamma}+
	R^\alpha_{\beta\delta\smbar\gamma}\,\psi^\delta\wedge\bar{\psi^\gamma},
	$$
where
$$\eqalign{R^\alpha_{\beta;\mu\smbar\nu}&=-\delta_{\smbar\nu}(\Gamma^\alpha
	_{\beta;\mu})-\Gamma^\alpha_{\beta\sigma}\delta_{\smbar\nu}(\Gamma^
	\sigma_{;\mu}),\cr
	R^\alpha_{\beta\delta;\smbar\nu}&=-\delta_{\smbar\nu}(\Gamma^\alpha
	_{\beta\delta})=R^\alpha_{\delta\beta;\smb\nu},\cr
	R^\alpha_{\beta\smbar\gamma;\mu}&=-\dot\de_{\smb\gamma}(\Gamma^\alpha
	_{\beta;\mu})-\Gamma^\alpha_{\beta\sigma}\Gamma^\sigma_{\smbar\gamma;
	\mu},\cr
	R^\alpha_{\beta\delta\smbar\gamma}&=-\dot\de_{\smb\gamma}(\Gamma^\alpha
	_{\beta\delta})=R^\alpha_{\delta\beta\smb\gamma}.\cr}\neweq\eqR$$
In particular, since
$$\eqalign{(D\tau)^\alpha=(\eta^V\wedge\Omega)^\alpha&=\psi^\sigma\wedge\Omega
	^\alpha_\sigma=
	R^\alpha_{\sigma;\mu\smbar\nu}\,\psi^\sigma\wedge dz^\mu\wedge d\smb z^
	\nu+R^\alpha_{\sigma\delta;\smbar\nu}\,\psi^\sigma\wedge\psi^\delta
	\wedge d\smb z^\nu\cr
	&\quad+R^\alpha_{\sigma\smbar\gamma;\mu}\,\psi^\sigma\wedge dz^\mu
	\wedge\bar{\psi^\gamma}+R^\alpha_{\sigma\delta\smb\gamma}\,
	\psi^\sigma\wedge\psi^\delta\wedge\bar{\psi^\gamma}\cr
	&=-R^\alpha_{\sigma;\mu\smbar\nu}\,dz^\mu\wedge\psi^\sigma\wedge 
	d\smb z^\nu-R^\alpha_{\sigma\smb\gamma;\mu}\,dz^\mu\wedge\psi^\sigma
	\wedge\bar{\psi^\gamma},\cr}$$
the vanishing of $\tau$ implies the vanishing of most of the curvature.
 
Another consequence of \eqR\ is an unexpected relation between $\Omega$ and 
$\tau$:
 
\newthm Lemma \luno: Let $D$ be the Chern-Finsler connection associated 
to a strongly pseudoconvex Finsler metric~$F$ on a complex 
manifold~$M$. Then
$$\tau=\Omega(\cdot,\,\cdot)\iota.$$
 
\pf Recalling \eqR, \eqdf, \eqde\ and  
$$\Gamma^\alpha_{\beta;\mu}v^\beta=\Gamma^\alpha_{;\mu},\qquad\Gamma^\alpha_
	{\beta\gamma}v^\beta=0,$$
we have
$$\displaylines{R^\alpha_{\beta;\mu\smbar\nu}v^\beta=-\delta_{\smbar\nu}
	(\Gamma^\alpha_{;\mu}),\cr
	R^\alpha_{\beta\delta;\smbar\nu}v^\beta=0,\cr
	R^\alpha_{\beta\smbar\gamma;\mu}v^\beta=-\Gamma^\alpha_{\smb\gamma;
	\mu},\cr
	R^\alpha_{\beta\delta\smbar\gamma}v^\beta=0,\cr}$$
and the assertion follows from \eqTau.\qedn
 
\smallsect 5. Holomorphic curvature
 
One of the most useful concept in hermitian geometry is the notion 
of holomorphic sectional curvature. To find the correct analogue in our 
setting, we first need a closer look to the horizontal part of the curvature 
operator. We define the {\sl horizontal curvature tensor}~$R$ by 
$$R_v(H,\bar K,L,\bar M)=\langle\Omega(H,\bar K)L,M\rangle_v$$
for all $H$, $K$, $L$, $M\in\ca H_v$ and $v\in\tildeM$. In local coordinates,
$$R(H,\bar K,L,\bar M)=G_{\sigma\smb\beta}R^\sigma_{\alpha;\mu\smb\nu}
	H^\mu\bar{K^\nu}L^\alpha\bar{M^\beta}.$$
The symmetries of~$R$ are easily described:
 
\newthm Proposition \tnove: Take $v\in\tildeM$ and $H$,~$K$,~$L$,~$M\in\ca 
H_v$. Then
$$R(\bar{K},H,L,\bar{M})=-R(H,\bar{K},L,\bar{M});\neweq\eqtquat$$
$$R(K,\bar{H},M,\bar{L})=\bar{R(H,\bar{K},L,\bar{M})}.\neweq\eqtcin$$
Furthermore, if $\smb\de_H\theta=0$ we also have
$$R(L,\bar{K},H,\bar{M})=R(H,\bar{K},L,\bar{M})=R(H,\bar{M},L,\bar{K}).
	\neweq\eqtsei$$
 
\pf \eqtquat\ follows immediately from the observation $\Omega^\alpha_\beta(
\bar{K},H)=-\Omega^\alpha_\beta(H,\bar{K})$. To prove \eqtcin, we start from
$$\Omega^\alpha_\beta=\smb\de\omega^\alpha_\beta=\smb\de(G^{\smbar\tau\alpha}
	\de G_{\beta\smbar\tau})=-G^{\smbar\tau\mu}G^{\smbar\nu\alpha}\smb\de
	G_{\mu\smbar\nu}\wedge\de G_{\beta\smbar\tau}+G^{\smbar\tau\alpha}\bar
	\de\de G_{\beta\smbar\tau};$$
in particular,
$$G_{\alpha\smbar\gamma}\Omega^\alpha_\beta=-G^{\smbar\tau\mu}\smb\de G_{\mu
	\smbar\gamma}\wedge\de G_{\beta\smbar\tau}+\smb\de\de G_{\beta\smbar
	\gamma}.$$
On the other hand,
$$\Omega^{\smbar\alpha}_{\smbar\gamma}=G^{\smbar\nu\tau}G^{\smbar\alpha\mu}\bar
	\de G_{\tau\smbar\gamma}\wedge\de G_{\mu\smbar\nu}-G^{\smbar\alpha\tau}
	\smb\de\de G_{\tau\smbar\gamma};$$
hence
$$G_{\alpha\smbar\gamma}\Omega^\alpha_\beta=-G_{\beta\smbar\alpha}
		\Omega^{\smbar\alpha}_{\smbar\gamma}.$$
In our case, this means that
$$\eqalign{R(K,\bar{H},M,\bar L)&=G_{\alpha\smbar\gamma}\Omega^\alpha_\beta(K,
	\bar H)M^\beta\bar{L^\gamma}=-G_{\beta\smb\alpha}\Omega^{\smb\alpha}
	_{\smb\gamma}(K,\bar{H})M^\beta\bar{L^\gamma}\cr
	&=\bar{G_{\alpha\smb\beta}\Omega^\alpha_\gamma(H,\bar K)L^\gamma
	\bar{M^\beta}}\cr
	&=\bar{R(H,\bar K,L,\bar M)},\cr}$$
and \eqtcin\ is proved.
 
Now, \eqtsei. First of all, $\smb\de_H\theta=p^*_H(D\theta)$, because we saw 
that the (2,0)-part of $D\theta$ vanishes. Proposition~\dB\ 
says that $D\theta=\eta^H\wedge\Omega$; in local coordinates, 
$$\eqalign{(\eta^H\wedge\Omega)^\alpha=dz^\sigma\wedge\Omega^\alpha_\sigma=&
	R^\alpha_{\sigma;\mu\smb\nu}\,dz^\sigma\wedge dz^\mu\wedge d\smb z^\nu+
	R^\alpha_{\sigma\delta;\smb\nu}\,dz^\sigma\wedge\psi^\delta\wedge 
	d\smb z^\nu\cr
	&+R^\alpha_{\sigma\smb\gamma;\mu}\,dz^\sigma\wedge dz^\mu\wedge\bar
	{\psi^\gamma}+R^\alpha_{\sigma\delta\smb\gamma}\,dz^\sigma\wedge
	\psi^\delta\wedge\bar{\psi^\gamma};\cr}\neweq\eqtth$$
in particular, $\smb\de_H\theta=0$ iff $R^\alpha_{\sigma;\mu\smb\nu}=R^
\alpha_{\mu;\sigma\smb\nu}$. Then
$$R(L,\bar K,H,\bar M)=G_{\alpha\smb\tau}R^\alpha_{\sigma;\mu\smb\nu}L^\mu
	\bar{K^\nu}H^\sigma\bar{M^\tau}=G_{\alpha\smb\tau}R^\alpha_{\mu;\sigma
	\smb\nu}L^\mu\bar{K^\nu}H^\sigma\bar{M^\tau}=R(H,\bar K,L,\bar M).$$
Finally,
$$R(H,\bar M,L,\bar K)=\bar{R(M,\bar H,K,\bar L)}=\bar{R(K,\bar H,M,\bar L)}
	=R(H,\bar K,L,\bar M).$$
\qedn
 
We remark that \eqtcin\ is equivalent to
$$\langle\Omega(H,\bar K)L,M\rangle=\langle L,\Omega(K,\bar H)M\rangle
	\neweq\eqtcb$$
for all~$H$, $K$, $L$, $M\in\ca H$.
 
Now, one possible approach to the holomorphic sectional curvature is to
consider the {\sl (horizontal) holomorphic flag curvature} $\tilde K_F(H)$
of~$F$ along a horizontal vector~$H\in \ca H_v$: 
$$\tilde K_F(H)={2\over\langle H,H\rangle^2_v}R(H,\bar H,H,\bar H).$$
Exactly as in the hermitian case, if $\smb\de_H\theta=0$ then the holomorphic 
flag curvature completely determines the horizontal curvature tensor:
 
\newthm Proposition \tdieci: Let $R$,~$S\colon\ca H_v\times\bar{\ca H_v}\times
\ca H_v\times\bar{\ca H_v}\to\C$ be two quadrilinear maps satisfying $\eqtcin$ 
and $\eqtsei$. Assume that
$$\forevery{H\in\ca H_v}R(H,\bar H,H,\bar H)=S(H,\bar H,H,\bar H).$$
Then $R\equiv S$.
 
The proof is very similar to the traditional one for hermitian metrics; 
see~[KN] and~[AP4] for the details. We do not discuss it here because, from a 
certain 
point of view, the holomorphic flag curvature is {\it not} the right 
generalization of the holomorphic sectional curvature. In fact, 
roughly speaking, it contains too many informations. Requiring, for instance, 
that the holomorphic flag curvature is constant means imposing very strong 
constraints on the behavior of the complex Finsler metric, constraints that 
are somewhat beyond the geometry of the metric which lives naturally on the 
tangent bundle of the manifold. Of course, one may study 
such requirements, but in this case the theory seems to be a standard 
consequence of the hermitian geometry of vector bundles without significant 
application to the function theory of the manifold.
 
A different notion appears to be a more appropriate tool for the applications
in complex geometry (see~[K3], [AP1], [AP2], [AP3], and sections 6 and 7 where 
we discuss variational formulas; cf. also~[Ch] and~[BC] for similar arguments 
in the real case). 
Namely, let $F\colon T^{1,0}M\to\R^+$ be a strongly pseudoconvex Finsler 
metric on a complex manifold~$M$, and take $v\in\tildeM$. 
Then the {\sl holomorphic curvature} $K_F(v)$ of $F$ along~$v$ is given by
$$K_F(v)=\tilde K_F\bigl(\chi(v)\bigr)={2\over G(v)^2}R\bigl(\chi(v),\smbar
	\chi(v),\chi(v),\smbar\chi(v)\bigr).$$
Clearly,
$$K_F(\zeta v)=K_F(v)$$
for all $\zeta\in\C^*$; so this is the holomorphic curvature discussed by 
Kobayashi~[K3]. Note that, by Proposition~\tnove, the holomorphic 
curvature is necessarily real-valued.
 
In local coordinates we get
$$K_F=-{2\over G^2}G_\alpha\delta_{\smb\nu}(\Gamma^\alpha_{;\mu})v^\mu
	\bar{v^\nu}.\neweq\eqKf$$
If $F$ comes from a hermitian metric, \eqKf\ gives exactly the classical
holomorphic sectional curvature. Furthermore, our definition recovers another
important geometrical characterization of the holomorphic sectional curvature,
and provides a firm link with the theory of invariant metrics on complex
manifolds (cf.~[AP2]). Wu~[Wu] has shown that for a hermitian metric~$g$ on a
complex manifold~$M$, the holomorphic sectional curvature of~$g$ along $v\in
T^{1,0}_p M$ is the maximum value attained by the Gaussian curvature at the
origin of the pull-back metric $\phe^*g$ when $\phe$ varies among the
holomorphic maps from the unit disk~$\Delta\subset\C$ into~$M$ with
$\phe(0)=p$ and $\phe'(0)= \lambda v$ for some $\lambda\in\C^*$. Well, this is
true in our case too: 
 
\newthm Theorem \holc: Let $F\colon T^{1,0}M\to\R^+$ be a strongly 
pseudoconvex Finsler metric on a complex manifold~$M$, and take $p\in M$ and 
$v\in\tildeM_p$. Then
$$K_F(v)=\sup\{K(\phe^*G)(0)\},$$
where $K(\phe^*G)(0)$ is the Gaussian curvature at the origin of the pull-back 
metric $\phe^*G$, and the supremum is taken with respect to the family of all 
holomorphic maps $\phe\colon\Delta\to M$ with $\phe(0)=p$ and $\phe'(0)=
\lambda v$ for some $\lambda\in\C^*$.
 
For the proof, see~[AP2]. We also recall that this variational interpretation 
of the holomorphic curvature makes sense for upper semicontinuous Finsler 
metrics, and has been previously investigated in geometric function theory 
(see~[W], [R] and~[S]).
 
This ends the general presentation of the setting we 
suggest for studying complex Finsler geometry. To substantiate this 
suggestion, in the next sections we give a few applications: the variation 
formulas and a close look to manifold with constant holomorphic curvature.
 
\smallsect 6. First variation of the length integral and geodesics
 
Let $F\colon T^{1,0}M\to\R^+$ be a strongly pseudoconvex Finsler 
metric on a complex manifold~$M$. To~$F$ we may associate a function $F^o\colon
T_\R M\to\R^+$ just by setting
$$\forevery{u\in T_\R M}F^o(u)=F(u_o),$$
where $u\mapsto u_o=(u-iJu)/2$ is the standard isomorphism between $T_\R M$ 
and $T^{1,0}M$ ($J$~is the complex structure on~$T_\R M$).
$F^o$ satisfies all the properties defining a real Finsler metric, but perhaps 
the indicatrices are not necessarily strongly convex. 
Nevertheless, we may use it to measure the length of curves, and so to define 
geodesics; and one of the main results of this section is a theorem ensuring 
the local existence and uniqueness of geodesics for weakly K\"ahler Finsler 
metrics only under the strong pseudoconvexity hypothesis --- a striking 
by-product of the complex structure.
 
Let us fix the notations needed to study variations of the length integral in 
this setting. The idea is, as usual, to pull back the connection along a 
curve; but since our connection lives on the tangent-tangent bundle, the 
details are a bit delicate.
 
A {\sl regular curve} $\sigma\colon[a,b]\to M$ is a $C^1$ curve with never 
vanishing tangent vector. Here, we mean the tangent vector in $T^{1,0}M$, 
obtained via the canonical isomorphism with $T_\R M$: so we set
$$\dot\sigma(t)=d^{1,0}\sigma_t\left({d\over dt}\right)={d\sigma^\alpha\over
	dt}(t)\left.{\de\over\de z^\alpha}\right|_{\sigma(t)},$$
where $d^{1,0}$ is the composition of the differential with the projection
of $T_\C M$ onto $T^{1,0}M$ associated to the splitting $T_\C M=T^{1,0}M
\oplus T^{0,1}M$.  
 
The {\sl length} of a regular curve~$\sigma$ with respect to the strongly
pseudoconvex Finsler metric~$F$ is given by
$$L(\sigma)=\int^b_a F\bigl(\dot\sigma(t)\bigr)\,dt,$$
exactly as in the hermitian case.
 
A geodesic for $F$ is a curve which is a critical point of the length 
functional. To be more precise, let $\sigma_0\colon[a,b]\to M$ be a regular 
curve with $F(\dot\sigma_0)\equiv c_0>0$. A {\sl regular variation} of~$\sigma
_0$ is a $C^1$ map $\Sigma\colon(-\eps,\eps)\times[a,b]\to M$ such that
\smallskip
\item{(a)} $\sigma_0(t)=\Sigma(0,t)$ for all $t\in[a,b]$;
\item{(b)} for every $s\in(-\eps,\eps)$ the curve $\sigma_s(t)=\Sigma(s,t)$ is 
a regular curve in~$M$;
\item{(c)} $F(\dot\sigma_s)\equiv c_s>0$ for every $s\in(-\eps,\eps)$.
\smallskip
\noindent A regular variation $\Sigma$ is {\sl fixed} if it moreover satisfies
\smallskip
\item{(d)} $\sigma_s(a)=\sigma_0(a)$ and $\sigma_s(b)=\sigma_0(b)$ for every
$s\in(-\eps,\eps)$.
\smallskip
\noindent If $\Sigma$ is a regular variation of the curve $\sigma_0$, we 
define the function $\ell_\Sigma\colon(-\eps,\eps)\to\R^+$ by
$$\ell_\Sigma(s)=L(\sigma_s).$$
We shall say that a regular curve $\sigma_0$ is a {\sl geodesic} for~$F$ iff
$${d\ell_\Sigma\over ds}(0)=0$$
for all fixed regular variations~$\Sigma$ of~$\sigma_0$.
 
Our first goal is to write the first variation of the length functional; we 
shall then find the differential equation satisfied by the geodesics (see 
also~[AP1]).
 
Let $\Sigma\colon(-\eps,\eps)\times[a,b]\to M$ be a regular variation of a 
regular curve $\sigma_0\colon[a,b]\to M$. Let $p\colon\Sigma^*(T^{1,0}M)\to
(-\eps,\eps)\times[a,b]$ be the 
pull-back bundle, and $\gamma\colon\Sigma^*(T^{1,0}M)\to T^{1,0}M$ be the 
bundle map such that the diagram
$$\dmatrix{\Sigma^*(T^{1,0}M)&\mapright{\gamma}&T^{1,0}M\cr
	\mapdownl{p}&&\mapdownr{\pi}\cr
	(-\eps,\eps)\times[a,b]&\mapright{\Sigma}&M\cr}$$
commutes. 
 
Two particularly important sections of $\Sigma^*(T^{1,0}M)$ are
$$T=\gamma^{-1}\left(d^{1,0}\Sigma\left({\de\over\de t}\right)\right)
	={\de\Sigma^\alpha\over\de t}{\de\over\de z^\alpha},$$
and
$$U=\gamma^{-1}\left(d^{1,0}\Sigma\left({\de\over\de s}\right)\right)
	={\de\Sigma^\alpha\over\de s}{\de\over\de z^\alpha};$$
the restriction of $U$ to $s=0$ is the {\sl variation vector} of the 
variation~$\Sigma$. Note that setting $\Sigma^*\tildeM=\gamma^{-1}(\tildeM)$, 
we have $T\in\ca X(\Sigma^*\tildeM)$ and
$$T(s,t)=\gamma^{-1}\bigl(\dot\sigma_s(t)\bigr).$$
Now we pull-back $T^{1,0}\tildeM$ over $\Sigma^*\tildeM$ by using~$\gamma$, 
obtaining the commutative diagram
$$\dmatrix{\gamma^*(T^{1,0}\tildeM)&\mapright{\tilde\gamma}&T^{1,0}\tildeM\cr
	\mapdownl{}&&\mapdownr{}\cr
	\Sigma^*\tildeM&\mapright{\gamma}&\tildeM\cr
	\mapdownl{}&&\mapdownr{}\cr
	(-\eps,\eps)\times[a,b]&\mapright{\Sigma}&M\cr};$$
note that $\gamma^*(T^{1,0}\tildeM)$ is a complex vector bundle over a real 
manifold. The bundle map~$\tilde\gamma$
induces a hermitian structure on $\gamma^*(T^{1,0}\tildeM)$ by 
$$\forevery{X,Y\in\gamma^*(T^{1,0}\tildeM)_v}\langle X,Y\rangle_v=\langle
	\tilde\gamma(X),\tilde\gamma(Y)\rangle_{\gamma(v)}.$$
Analogously, the Chern connection $D$ gives rise to a (1,0)-connection 
$$D^*\colon\ca X\bigl(\gamma^*(T^{1,0}\tildeM)\bigr)\to\ca X\bigl(T^*_\C(
	\Sigma^*\tildeM)\otimes\gamma^*(T^{1,0}\tildeM)\bigr),$$ 
where $T^*_\C(\Sigma^*\tildeM)=T^*_\R(\Sigma^*\tildeM)\otimes\C$, by setting
$$\eqalign{\nabla^*_{\!X}Y&=\tilde\gamma^{-1}\left(\nabla_{d^{1,0}\gamma(X)}
	\tilde\gamma(Y)\right),\cr
	\nabla^*_{\bar X}Y&=\tilde\gamma^{-1}\left(\nabla_{\bar{d^{1,0}
	\gamma(X)}}\,\tilde\gamma(Y)\right),\cr}$$
for all $X\in T_\R(\Sigma^*\tildeM)$ and $Y\in\ca X\bigl(\gamma^*(T^{1,0}
\tildeM)\bigr)$. In particular we have
$$\eqalign{X\langle Y,Z\rangle&=X\bigl(\langle\tilde\gamma(Y),\tilde\gamma(Z)
	\rangle_\gamma\bigr)=d\gamma(X)\bigl(\langle\tilde\gamma(Y),
	\tilde\gamma(Z)\rangle\bigr)\cr
	&=(d^{1,0}\gamma(X)+\bar{d^{1,0}\gamma(X)})\bigl(\langle\tilde\gamma(Y),
	\tilde\gamma(Z)\rangle\bigr)\cr
	&=\langle\nabla^*_{\!X}Y,Z\rangle+\langle Y,\nabla^*_{\bar X}Z\rangle
	+\langle\nabla^*_{\bar X}Y,Z\rangle+\langle Y,\nabla^*_{\!X}Z\rangle,
	\cr}\neweq\eqquno$$
for all $X\in T_\R(\Sigma^*\tildeM)$ and $Y$,~$Z\in\ca X\bigl(\gamma^*(T^{1,0}
\tildeM)\bigr)$.
 
We may also decompose $T_\R(\Sigma^*\tildeM)=\ca H^*\oplus\ca V^*$, 
where as usual a local real frame for~$\ca V^*$ is given by $\{\dot\de_\alpha,i
\dot\de_\alpha\}$, and a local frame for $\ca H^*$ is given by
$$\delta_t=\de_t-(\Gamma^\mu_{;\alpha}\circ\gamma){\de\Sigma^\alpha\over\de t}
	\dot\de_\mu,\qquad\delta_s=\de_s-(\Gamma^\mu_{;\alpha}\circ\gamma)
	{\de\Sigma^\alpha\over\de s}\dot\de_\mu,$$
where $\de_t=\de/\de t$ and $\de_s=\de/\de s$.
Therefore, setting $T^H=d^{1,0}\gamma(\delta_t)$ and $U^H=d^{1,0}\gamma
(\delta_s)$, we have
$$T^H(v)={\de\Sigma^\mu\over\de t}(s,t)\delta_\mu|_{\gamma(v)}
	=\chi_{\gamma(v)}\bigl(\dot\sigma_s(t)\bigr)\in\ca H_{\gamma(v)}$$
and
$$U^H(v)={\de\Sigma^\mu\over\de s}(s,t)\delta_\mu|_{\gamma(v)}
	=\chi_{\gamma(v)}\bigl(\gamma\bigl(U(s,t)\bigr)\bigr)\in\ca 
	H_{\gamma(v)},$$
for all $v\in\Sigma^*\tildeM_{(s,t)}$; they are the horizontal lifts of 
$\gamma(T)$ and $\gamma(U)$ respectively. In particular,
$$T^H\bigl(\gamma^{-1}(\dot\sigma_s)\bigr)=\chi(\dot\sigma_s).\neweq\eqqunom$$
 
If we take $v\in(\Sigma^*\tildeM)_{(s,t)}$, then 
$$d^{1,0}\gamma_v\bigl(T_\R(\Sigma^*M)\bigr)\subset T^{1,0}_{\gamma(v)}\tildeM
	\qquad\hbox{and}\qquad\tilde\gamma\bigl(\gamma^*(T^{1,0}\tildeM)_v
	\bigr)=T^{1,0}_{\tilde\gamma(v)}\tildeM.$$
Therefore
we also have a bundle map $\Xi\colon T_\R(\Sigma^*\tildeM)\to\gamma^*(T^{1,0}
\tildeM)$ such that the diagram
$$\dmatrix{T_\R(\Sigma^*\tildeM)&\mapright{\Xi}&\gamma^*(T^{1,0}\tildeM)\cr
	&\mapse{d^{1,0}\gamma}&\mapdownr{\tilde\gamma}\cr
	&&T^{1,0}\tildeM\cr}$$
commutes. Using $\Xi$ we may prove three final formulas:
$$\eqalign{\tilde\gamma\bigl(\nabla^*_{\!X}\Xi(Y)-\nabla^*_{\!Y}\Xi(X)\bigr)&=
	\nablab_{d^{1,0}\gamma(X)}\,d^{1,0}\gamma(Y)-\nablab_{d^{1,0}\gamma(Y)}
	\,d^{1,0}\gamma(X)\cr
	&=\left[d^{1,0}\gamma(X),d^{1,0}\gamma(Y)\right]+\theta\bigl(d^{1,0}
	\gamma(X),d^{1,0}\gamma(Y)\bigr),\cr}\neweq\eqqdue$$
for all $X$,~$Y\in\ca X\bigl(T_\R(\Sigma^*\tildeM)\bigr)$;
$$\eqalign{\tilde\gamma\circ(\nabla^*_{\!X}\nabla^*_{\!Y}-\nabla^*_{\!Y}
	\nabla^*_{\!X})&=(\nablab_{d^{1,0}\gamma(X)}\nablab_{d^{1,0}\gamma(Y)}-
	\nablab_{d^{1,0}\gamma(Y)}\nablab_{d^{1,0}\gamma(X)})\circ\tilde
	\gamma\cr
	&=\nabla_{[d^{1,0}\gamma(X),d^{1,0}\gamma(Y)]}\circ\tilde\gamma,\cr}
	\neweq\eqqtre$$
and
$$\tilde\gamma\circ(\nabla^*_{\!X}\nabla^*_{\bar Y}-\nabla^*_{\bar Y}\nabla^*
	_{\!X})=\bigl(\nabla^*_{[d^{1,0}\gamma(X),\bar{d^{1,0}\gamma(Y)}]}+
	\Omega\bigl(d^{1,0}\gamma(X),\bar{d^{1,0}\gamma(Y)}\bigr)\circ\tilde
	\gamma,\neweq\eqqqua$$
for all $X$,~$Y\in T_\R(\Sigma^*\tildeM)$.
 
We are now able to prove the first variation formula for weakly
K\"ahler Finsler metrics:
 
\newthm Theorem \FVFK: Let $F\colon T^{1,0}M\to\R^+$ be a weakly 
K\"ahler Finsler metric on a complex manifold~$M$. Take a regular curve 
$\sigma_0\colon[a,b]\to M$ with $F(\dot\sigma_0)\equiv c_0>0$, and a regular 
variation $\Sigma\colon(-\eps,\eps)\times[a,b]\to M$ of~$\sigma_0$. Then
$${d\ell_\Sigma\over ds}(0)={1\over c_0}\left\{\Re\langle U^H,T^H\rangle
	_{\dot\sigma_0}\biggr|^b_a-\Re\int^b_a\langle U^H,\nablab_{T^H+
	\bar{T^H}}T^H\rangle_{\dot\sigma_0}\,dt\right\}.$$
In particular, if $\Sigma$ is a fixed variation, that is $\Sigma(\cdot,a)\equiv
\sigma_0(a)$ and $\Sigma(\cdot,b)\equiv\sigma_0(b)$, we have
$${d\ell_\Sigma\over ds}(0)=-{1\over c_0}\Re\int^b_a\langle U^H,
	\nablab_{T^H+\bar{T^H}}T^H\rangle_{\dot\sigma_0}\,dt.\neweq\eqFixv$$
 
\pf By definition,
$$\ell_\Sigma(s)=\int^b_a\bigl(G(\dot\sigma_s)\bigr)^{1/2}\,dt;$$
therefore
$${d\ell_\Sigma\over ds}={1\over 2c_s}\int^b_a{\de\over\de s}[G(\dot\sigma_s)]
	\,dt={1\over 2c_s}\int^b_a{\de\over\de s}\langle\Xi(\delta_t),
	\Xi(\delta_t)\rangle_T\,dt,$$
where $c_s\equiv F(\dot\sigma_s)$ and we used
$$G(\dot\sigma_s)=\langle\chi(\dot\sigma_s),\chi(\dot\sigma_s)\rangle_{\dot
	\sigma_s}=\langle\Xi(\delta_t),\Xi(\delta_t)\rangle_T,$$
by \eqqunom. Now, using \eqquno\ and \eqqdue, we get
$$\eqalign{{1\over2}&{\de\over\de s}\langle\Xi(\delta_t),\Xi(\delta_t)\rangle_T
	={1\over2}\delta_s\langle\Xi(\delta_t),\Xi(\delta_t)\rangle_T\cr
	&={1\over2}\biggl\{\langle\nabla^*_{\delta_s}\Xi(\delta_t),
	\Xi(\delta_t)\rangle_T+\langle\Xi(\delta_t),\nabla^*_{\bar{\delta_s}}
	\,\Xi(\delta_t)\rangle_T\cr
	&\qquad\phantom{{1\over2}\biggl\{\langle\nabla^*_{\delta_s}\Xi(\delta_t)
	,\Xi(\delta_t)\rangle_T}+\langle\nabla^*_{\bar{\delta_s}}\,
	\Xi(\delta_t),\Xi(\delta_t)\rangle_T+\langle\Xi(\delta_t),\nabla^*
	_{\delta_s}\Xi(\delta_t)\rangle_T\biggr\}\cr
	&=\Re\left\{\langle\nabla^*_{\delta_s}\Xi(\delta_t),\Xi(\delta_t)
	\rangle_T+\langle\nabla^*_{\bar{\delta_s}}\,\Xi(\delta_t),\Xi(\delta_t)
	\rangle_T\right\}\cr
	&=\Re\left\{\langle\nabla^*_{\delta_t}\Xi(\delta_s),\Xi(\delta_t)
	\rangle_T+\langle[U^H,T^H]+\nabla_{\bar{U^H}}T^H,T^H\rangle_{\dot
	\sigma_s}+\langle\theta(U^H,T^H),T^H\rangle_{\dot\sigma_s}\right\}.
	\cr}$$
Since $F$ is weakly K\"ahler, \eqqunom\ yields
$$\langle \theta(U^H,T^H),T^H\rangle_{\dot\sigma_s}=0.$$
Furthermore,
$$\eqalign{[U^H,T^H]&=\left\{{\de\Sigma^\nu\over\de s}\,\delta_\nu\!\left(
	{\de\Sigma^\mu\over\de t}\right)-{\de\Sigma^\nu\over\de t}\,\delta_\nu
	\!\left({\de\Sigma^\mu\over\de s}\right)\right\}\,\delta_\mu,\cr
	\nabla_{\bar{U^H}}T^H&={\de\bar{\Sigma^\nu}\over\de s}\,\delta_{\smb\nu}
	\!\left({\de\Sigma^\mu\over\de t}\right)\delta_\mu;\cr}$$
since
$${\de\Sigma^\nu\over\de s}\,\delta_\nu\!\left({\de\Sigma^\mu
	\over\de t}\right)+{\de\bar{\Sigma^\nu}\over\de s}\,\delta_{\smb\nu}\!
	\left({\de\Sigma^\mu\over\de t}\right)={\de^2\Sigma^\mu\over\de s\de t}
	={\de^2\Sigma^\mu\over\de t\de s}={\de\Sigma^\nu\over\de t}\,\delta_\nu
	\!\left({\de\Sigma^\mu\over\de s}\right)+{\de\bar{\Sigma^\nu}\over\de t}
	\,\delta_{\smb\nu}\!\left({\de\Sigma^\mu\over\de s}\right),
	\neweq\eqts$$
we get
$$[U^H,T^H]+\nabla_{\bar{U^H}}T^H={\de\bar{\Sigma^\nu}\over\de t}\,\delta_{\smb
	\nu}\!\left({\de\Sigma^\mu\over\de s}\right)\delta_\mu=\nabla
	_{\bar{T^H}}U^H.\neweq\eqqcin$$
Then
$$\eqalign{{1\over2}{\de\over\de s}&\langle\Xi(\delta_t),\Xi(\delta_t)\rangle_T
	=\Re\left\{\langle\nabla^*_{\delta_t}\Xi(\delta_s),\Xi(\delta_t)
	\rangle_T+\langle\nabla^*_{\bar{\delta_t}}\,\Xi(\delta_s),\Xi(\delta_t)
	\rangle_T\right\}\cr
	&=\Re\left\{\delta_t\langle\Xi(\delta_s),\Xi(\delta_t)\rangle_T-
	\langle\Xi(\delta_s),\nabla^*_{\!\delta_t+\bar{\delta_t}}\,\Xi(\delta_t)
	\rangle_T\right\}\cr
	&=\Re\left\{{\de\over\de t}\langle U^H,T^H\rangle_{\dot\sigma_s}-
	\langle U^H,\nablab_{T^H+\bar{T^H}}T^H\rangle_{\dot\sigma_s}\right\},
	\cr}\neweq\eqcp$$
and the assertion follows.\qedn
 
As a corollary we get the equation of geodesics:
 
\newthm Corollary \geod: Let $F\colon T^{1,0}M\to\R^+$ be a weakly 
K\"ahler Finsler metric on a complex manifold~$M$, and let $\sigma\colon[a,b]
\to M$ be a regular curve with~$F(\dot\sigma)\equiv c_0>0$. Then $\sigma$~is a 
geodesic for~$F$ iff
$$\nablab_{T^H+\bar{T^H}}T^H\equiv0,\neweq\eqgeod$$
where $T^H(v)=\chi_v\bigl(\dot\sigma(t)\bigr)\in\ca H_v$ for all $v\in\tildeM_
{\sigma(t)}$.
 
\pf It follows immediately from \eqFixv.\qedn
 
\newthm Corollary \geodb: Let $F\colon T^{1,0}M\to\R^+$ be a weakly 
K\"ahler Finsler metric on a complex manifold~$M$. Then for any $p\in M$ and 
$v\in\tildeM_p$ with $F(v)=1$ there exists a unique geodesic $\sigma\colon
(-\eps,\eps)\to M$ such that $\sigma(0)=p$ and $\dot\sigma(0)=v$.
 
\pf In local coordinates we have
$$\nablab_{T^H+\bar{T^H}}T^H=\left[(\dot\sigma^\mu\delta_\mu+\bar{\dot
	\sigma^\mu}
	\delta_{\smb\mu})(\dot\sigma^\alpha)+\Gamma^\alpha_{\nu;\mu}(\dot\sigma)
	\dot\sigma^\mu\dot\sigma^\nu\right]\delta_\alpha=[\ddot\sigma^\alpha
	+\Gamma^\alpha_{;\mu}(\dot\sigma)\dot\sigma^\mu]\delta_\alpha.$$
So \eqgeod\ is a quasi-linear O.D.E. system, and the assertion follows.\qedn
 
Thus the standard O.D.E. arguments apply in this case too, and we may 
recover for weakly K\"ahler Finsler metrics the usual theory of geodesics.
In particular, if the metric~$F$ is complete we can define the {\sl 
exponential map} $\exp_p\colon T^{1,0}_pM\to M$ for any~$p\in M$. See~[AP4] 
for details.
 
\smallsect 7. Second variation of the length integral
 
Our next goal is the second variation formula, which holds for 
K\"ahler Finsler metrics. To express it correctly, we need two further 
ingredients. The first one is the {\sl horizontal (1,1)-torsion} $\tau^H$, 
simply defined by
$$\tau^H(X,\bar Y)=\Theta\bigl(\tau(X,\bar Y)\bigr)=\Omega(X,\bar Y)\chi.$$
The second one is the {\sl symmetric product} $\llangle\,,\rrangle\colon
\ca H\times\ca H\to\C$ locally given by
$$\forevery{H,K\in\ca H_v}\llangle H,K\rrangle_v=G_{\alpha\beta}(v)\,H^\alpha
	K^\beta.$$
It is clearly globally well-defined, and it satisfies
$$\forevery{H\in\ca H}\llangle H,\chi\rrangle=0.$$
 
\newthm Theorem \SVFK: Let $F\colon T^{1,0}M\to\R^+$ be a 
K\"ahler Finsler metric on a complex manifold~$M$. Take a geodesic
$\sigma_0\colon[a,b]\to M$ with $F(\dot\sigma_0)\equiv1$, and let
$\Sigma\colon(-\eps,\eps)\times[a,b]\to M$ be a regular variation of~$\sigma_
0$. Then
$$\eqalign{{d^2\ell_\Sigma\over ds^2}(0)&=\Re\langle\nablab_{U^H+\bar{U^H}}
	U^H,T^H\rangle_{\dot\sigma_0}\biggr|^b_a\cr
	&\quad+\int^b_a\biggl\{\bigl\|\nablab_{T^H+\bar{T^H}}U^H\bigr\|^2
	_{\dot\sigma_0}-\left|{\de\over\de t}\Re\langle U^H,T^H\rangle_{\dot
	\sigma_0}\right|^2\cr
	&\qquad\qquad-\Re\Bigl[\langle\Omega(T^H,\bar{U^H})U^H,T^H\rangle_{\dot
	\sigma_0}-\langle\Omega(U^H,\bar{T^H})U^H,T^H\rangle_{\dot\sigma_0}\cr
	&\quad\qquad\qquad\qquad+\llangle\tau^H(U^H,\bar{T^H}),U^H\rrangle_{\dot
	\sigma_0}-\llangle\tau^H(T^H,\bar{U^H}),U^H\rrangle_{\dot\sigma_0}\Bigr]
	\biggr\}\,dt.
	\cr}$$
In particular, if $\Sigma$ is a fixed variation such that $\Re\langle U^H,T^H
\rangle_{\dot\sigma_0}$ is constant we have
$$\eqalign{{d^2\ell_\Sigma\over ds^2}(0)&=\int^b_a\biggl\{\bigl\|\nablab_{T^H+
	\bar{T^H}}U^H\bigr\|^2_{\dot\sigma_0}\cr
	&\qquad\qquad-\Re\Bigl[\langle\Omega(T^H,\bar{U^H})U^H,T^H\rangle_{\dot
	\sigma_0}-\langle\Omega(U^H,\bar{T^H})U^H,T^H\rangle_{\dot\sigma_0}\cr
	&\quad\qquad\qquad\qquad+\llangle\tau^H(U^H,\bar{T^H}),U^H\rrangle_{\dot
	\sigma_0}-\llangle\tau^H(T^H,\bar{U^H}),U^H\rrangle_{\dot\sigma_0}\Bigr]
	\biggr\}\,dt.\cr}$$
 
\pf During the proof of the first variation formula --- in \eqcp\ --- we saw 
that
$${d\ell_\Sigma\over ds}(s)=\Re\int^b_a{\langle\nabla^*_{\delta_t}
	\Xi(\delta_s),\Xi(\delta_t)\rangle_T+\langle\nabla^*_{\bar{\delta_t}}
	\,\Xi(\delta_s),\Xi(\delta_t)\rangle_T\over\bigl(\langle\Xi(\delta_t),
	\Xi(\delta_t)\rangle_T\bigr)^{1/2}}\,dt.$$
So we need to compute
$$\eqalign{{\de\over\de s}\biggl[&{\langle\nabla^*_{\delta_t}
	\Xi(\delta_s),\Xi(\delta_t)\rangle_T+\langle\nabla^*_{\bar{\delta_t}}
	\,\Xi(\delta_s),\Xi(\delta_t)\rangle_T\over\bigl(\langle\Xi(\delta_t),
	\Xi(\delta_t)\rangle_T\bigr)^{1/2}}\biggr]\cr
	&\qquad={\delta_s\langle\nabla^*_{\delta_t}\Xi(\delta_s),\Xi(\delta_t)
	\rangle_T+\delta_s\langle\nabla^*_{\bar{\delta_t}}\,
	\Xi(\delta_s),\Xi(\delta_t)\rangle_T\over\bigl(\langle\Xi(\delta_t),
	\Xi(\delta_t)\rangle_T\bigr)^{1/2}}\cr
	&\qquad\qquad-{1\over2}{\langle\nabla^*
	_{\delta_t}\Xi(\delta_s),\Xi(\delta_t)\rangle_T+\langle\nabla^*
	_{\bar{\delta_t}}\,\Xi(\delta_s),\Xi(\delta_t)\rangle_T\over\bigl(
	\langle\Xi(\delta_t),\Xi(\delta_t)\rangle_T\bigr)^{3/2}}\,\delta_s
	\langle\Xi(\delta_t),\Xi(\delta_t)\rangle_T.\cr}\neweq\eqqfu$$
Since, when $s=0$, the denominator of the first term is equal to~1, and the 
denominator of the second term is equal to~2, we may forget them. Let us call 
(I) the numerator of the first term, and (I\negthinspace I) the numerator of 
the second term. First of all, \eqcp\ yields
$${1\over2}\Re\hbox{(I\negthinspace I)}=\left|\Re\left[{\de\over\de t}
	\langle U^H,T^H\rangle_{\dot\sigma_s}-\langle U^H,\nablab_{T^H+
	\bar{T^H}}T^H\rangle_{\dot\sigma_s}\right]\right|^2;$$
in particular, for $s=0$ we get
$${1\over2}\Re\hbox{(I\negthinspace I)}(0)=\left|{\de\over\de t}\Re
	\langle U^H,T^H\rangle_{\dot\sigma_0}\right|^2,\neweq\eqqfd$$
because $\sigma_0$ is a geodesic.
 
The computation of (I) is quite longer. First of all, using \eqqdue, \eqqtre\ 
and \eqqqua\ we get
$$\eqalign{\hbox{(I)}&=\langle\nabla^*_{\delta_s}\nabla^*_{\delta_t}
	\Xi(\delta_s),\Xi(\delta_t)\rangle_T+\langle\nabla^*_{\bar{\delta_s}}
	\nabla^*_{\delta_t}\Xi(\delta_s),\Xi(\delta_t)\rangle_T\cr
	&\quad+\langle\nabla^*_{\delta_t}\Xi(\delta_s),\nabla^*_{\delta_s}
	\Xi(\delta_t)\rangle_T+\langle\nabla^*_{\delta_t}\Xi(\delta_s),
	\nabla^*_{\bar{\delta_s}}\,\Xi(\delta_t)\rangle_T\cr
	&\quad+\langle\nabla^*_{\delta_s}\nabla^*_{\bar{\delta_t}}\,
	\Xi(\delta_s),\Xi(\delta_t)\rangle_T+\langle\nabla^*_{\bar{\delta_s}}
	\nabla^*_{\bar{\delta_t}}\,\Xi(\delta_s),\Xi(\delta_t)\rangle_T\cr
	&\quad+\langle\nabla^*_{\bar{\delta_t}}\,\Xi(\delta_s),
	\nabla^*_{\delta_s}\Xi(\delta_t)\rangle_T+\langle\nabla^*_{\bar
	{\delta_t}}\,\Xi(\delta_s),\nabla^*_{\bar{\delta_s}}\,\Xi(\delta_t)
	\rangle_T\cr
	&=\langle\nabla^*_{\delta_t}\nabla^*_{\delta_s}\Xi(\delta_s),
	\Xi(\delta_t)\rangle_T-\langle\nabla_{[T^H,U^H]}U^H,T^H\rangle_{\dot
	\sigma_s}\cr
	&\quad+\langle\nabla^*_{\delta_t}\nabla^*_{\bar{\delta_s}}\,
	\Xi(\delta_s),
	\Xi(\delta_t)\rangle_T-\langle\nabla_{[T^H,\bar{U^H}]}U^H,T^H\rangle
	_{\dot\sigma_s}-\langle\Omega(T^H,\bar{U^H})U^H,T^H\rangle_{\dot
	\sigma_s}\cr
	&\quad+\langle\nabla^*_{\delta_t}\Xi(\delta_s),\nabla^*_{\delta_t}
	\Xi(\delta_s)\rangle_T+\langle\nablab_{T^H}U^H,[U^H,T^H]+\nabla_{\bar
	{U^H}}T^H\rangle_{\dot\sigma_s}\cr
	&\quad\phantom{\quad+\langle\nabla^*_{\delta_t}\Xi(\delta_s),
	\nabla^*_{\delta_t}
	\Xi(\delta_s)\rangle_T}+\langle\nablab_{T^H}U^H,\theta(U^H,T^H)
	\rangle_{\dot\sigma_s}\cr
	&\quad+\langle\nabla^*_{\bar{\delta_t}}\nabla^*_{\delta_s}\Xi(\delta_s),
	\Xi(\delta_t)\rangle_T-\langle\nabla_{[\bar{T^H},U^H]}U^H,T^H\rangle
	_{\dot\sigma_s}+\langle\Omega(U^H,\bar{T^H})U^H,T^H\rangle_{\dot
	\sigma_s}\cr
	&\quad+\langle\nabla^*_{\bar{\delta_t}}\nabla^*_{\bar{\delta_s}}\,
	\Xi(\delta_s),\Xi(\delta_t)\rangle_T-\langle\nabla_{[\bar{T^H},
	\bar{U^H}]}U^H,T^H\rangle_{\dot\sigma_s}\cr
	&\quad+\langle\nabla^*_{\bar{\delta_t}}\,\Xi(\delta_s),\nabla^*
	_{\delta_t}\Xi(\delta_s)\rangle_T+\langle\nabla_{\bar{T^H}}U^H,[U^H,T^H]
	+\nabla_{\bar{U^H}}T^H\rangle_{\dot\sigma_s}\cr
	&\quad\phantom{\quad+\langle\nabla^*_{\bar{\delta_t}}\,\Xi(\delta_s),
	\nabla^*_{\delta_t}\Xi(\delta_s)\rangle_T}+\langle\nabla_{\bar{T^H}}
	U^H,\theta(U^H,T^H)\rangle_{\dot\sigma_s}.\cr}$$
Recalling \eqqcin, \eqqunom\ and that $F$ is K\"ahler we get
$$\eqalign{\hbox{(I)}&=\langle\nabla^*_{\delta_t}\nabla^*_{\delta_s}
	\Xi(\delta_s),\Xi(\delta_t)\rangle_T+\langle\nabla^*_{\delta_t}\nabla^*
	_{\bar{\delta_s}}\,\Xi(\delta_s),\Xi(\delta_t)\rangle_T\cr
	&\quad+\langle\nabla^*_{\bar{\delta_t}}\nabla^*_{\delta_s}
	\Xi(\delta_s),\Xi(\delta_t)\rangle_T+\langle\nabla^*_{\bar{\delta_t}}
	\nabla^*_{\bar{\delta_s}}\,\Xi(\delta_s),\Xi(\delta_t)\rangle_T\cr
	&\quad-\langle\Omega(T^H,\bar{U^H})U^H,T^H\rangle_{\dot\sigma_s}
	+\langle\Omega(U^H,\bar{T^H})U^H,T^H\rangle_{\dot\sigma_s}\cr
	&\quad-\langle\nabla_{[T^H,U^H]+[T^H,\bar{U^H}]+[\bar{T^H},U^H]+
	[\bar{T^H},\bar{U^H}]}U^H,T^H\rangle_{\dot\sigma_s}\cr
	&\quad+\langle\nabla^*_{\delta_t}\Xi(\delta_s),\nabla^*_{\delta_t}
	\Xi(\delta_s)\rangle_T+\langle\nabla^*_{\bar{\delta_t}}\,\Xi(\delta_s),
	\nabla^*_{\delta_t}\Xi(\delta_s)\rangle_T\cr
	&\quad+\langle\nabla^*_{\delta_t}\Xi(\delta_s),\nabla^*_{\bar{\delta_t}}
	\,\Xi(\delta_s)\rangle_T+\langle\nabla^*_{\bar{\delta_t}}\,
	\Xi(\delta_s),\nabla^*_{\bar{\delta_t}}\,\Xi(\delta_s)\rangle_T.\cr}$$
Now
$$\eqalign{[T^H,U^H]&=\left[{\de\Sigma^\mu\over\de t}\,\delta_\mu\!\left({\de
	\Sigma^\nu\over\de s}\right)-{\de\Sigma^\mu\over\de s}\,\delta_\mu\!
	\left({\de\Sigma^\nu\over\de t}\right)\right]\delta_\nu,\cr
	[T^H,\bar{U^H}]&={\de\Sigma^\mu\over\de t}{\de\bar{\Sigma^\nu}\over
	\de s}\bigl[\delta_{\smb\nu}(\Gamma^\alpha_{;\mu})\dot\de_\alpha
	-\delta_\mu(\Gamma^{\smb\beta}_{;\smb\nu})\dot\de_{\smb\beta}\bigr]
	+{\de\Sigma^\mu\over\de t}\,\delta_\mu\!\left({\de\bar{\Sigma^\nu}\over
	\de s}\right)\delta_{\smb\nu}-{\de\bar{\Sigma^\nu}\over\de s}\,
	\delta_{\smb\nu}\!\left({\de\Sigma^\mu\over\de t}\right)\delta_\mu,\cr
	[\bar{T^H},U^H]&={\de\bar{\Sigma^\mu}\over\de t}{\de\Sigma^\nu\over
	\de s}\bigl[\delta_\nu(\Gamma^{\smb\beta}_{;\smb\mu})\dot\de_{\smb\beta}
	-\delta_{\smb\mu}(\Gamma^\alpha_{;\nu})\dot\de_\alpha\bigr]
	+{\de\bar{\Sigma^\mu}\over\de t}\,\delta_{\smb\mu}\!\left({\de\Sigma^\nu
	\over\de s}\right)\delta_\nu-{\de\Sigma^\nu\over\de s}\,\delta_\nu\!
	\left({\de\bar{\Sigma^\mu}\over\de t}\right)\delta_{\smb\mu},\cr
	[\bar{T^H},\bar{U^H}]&=\left[{\de\bar{\Sigma^\mu}\over\de t}\,\delta
	_{\smb\mu}\!\left({\de\bar{\Sigma^\nu}\over\de s}\right)-{\de\bar{\Sigma
	^\mu}\over\de s}\,\delta_{\smb\mu}\!\left({\de\bar{\Sigma^\nu}\over\de 
	t}\right)\right]\delta_{\smb\nu},\cr}$$
and so \eqts\ and \eqTau\ yield
$$\displaylines{\qquad[T^H,U^H]+[T^H,\bar{U^H}]+[\bar{T^H},U^H]+[\bar{T^H},
	\bar{U^H}]\hfill\cr
	\hfill=\tau(U^H,\bar{T^H})-\tau(T^H,\bar{U^H})+\bar{\tau(U^H,\bar{T^H})}
	-\bar{\tau(T^H,\bar{U^H})}.\qquad\cr}$$
Furthermore, if $V\in\ca V$ we have
$$\eqalign{\langle\nablab_V U^H,T^H\rangle_{\dot\sigma_s}
	&=G_{\alpha\smb\beta}(\dot\sigma_s)V^\gamma\left[\dot\de_\gamma
	\!\left({\de\Sigma^\alpha\over\de s}\right)+\Gamma^\alpha
	_{\delta\gamma}(\dot\sigma_s)\left({\de\Sigma^\delta\over\de s}
	\right)\right]\bar{\dot\sigma^\beta_s}\cr
	&=G_\alpha(\dot\sigma_s)\Gamma^\alpha_{\delta\gamma}(\dot\sigma_s)
	\left({\de\Sigma^\delta\over\de s}\right)V^\gamma=
	G_{\delta\gamma}(\dot\sigma_s)\left({\de\Sigma^\delta\over\de s}
	\right)V^\gamma\cr
	&=\llangle\Theta(V),U^H\rrangle_{\dot\sigma_s},\cr}$$
and
$$\eqalign{\langle\nabla_{\bar V}U^H,T^H\rangle_{\dot\sigma_s}
	&=G_{\alpha\smb\beta}(\dot\sigma_s)\bar{V^\gamma}\left[\dot\de_
	{\smb\gamma}\!\left({\de\Sigma^\alpha\over\de s}\right)\right]
	\bar{\dot\sigma^\beta_s}\cr
	&=0.\cr}$$
Therefore
$$\displaylines{\qquad\langle\nabla_{[T^H,U^H]+[T^H,\bar{U^H}]+[\bar{T^H},U^H]+
	[\bar{T^H},\bar{U^H}]}U^H,T^H\rangle_{\dot\sigma_s}\hfill\cr
	\hfill=\llangle\tau^H(U^H,\bar{T^H}),U^H
	\rrangle_{\dot\sigma_s}-\llangle\tau^H(T^H,\bar{U^H}),U^H\rrangle
	_{\dot\sigma_s},\qquad\cr}$$
and thus
$$\eqalign{\hbox{(I)}&=\delta_t\langle\nabla^*_{\!\delta_s+\bar{\delta_s}}\,\Xi
	(\delta_s),\Xi(\delta_t)\rangle_T-\langle\nabla^*_{\!\delta_s+
	\bar{\delta_s}}\,\Xi(\delta_s),\nabla^*_{\!\delta_t+\bar{\delta_t}}
	\,\Xi(\delta_t)\rangle_T\cr
	&\quad-\langle\Omega(T^H,\bar{U^H})U^H,T^H\rangle_{\dot\sigma_s}
	+\langle\Omega(U^H,\bar{T^H}U^H,T^H\rangle_{\dot\sigma_s}\cr
	&\quad-\llangle\tau^H(U^H,\bar{T^H}),U^H\rrangle_{\dot\sigma_s}
	+\llangle\tau^H(T^H,\bar{U^H}),U^H\rrangle_{\dot\sigma_s}
	+\|\nabla^*_{\!\delta_t+\bar{\delta_t}}\,\Xi(\delta_s)\|^2_T.\cr}
	\neweq\eqqft$$
Recalling that for $s=0$ we have $\nabla^*_{\!\delta_t+\bar{\delta_t}}\,\Xi
(\delta_t)\equiv0$ because $\sigma_0$ is a geodesic, \eqqfu, \eqqfd\ and 
\eqqft\ evaluated at~$s=0$ yield the assertion.\qedn
 
So we have obtained the second variation formula for strongly pseudoconvex 
K\"ahler Finsler metrics. Besides its own intrinsic interest, we need it 
to compare the curvature of the real Finsler metric~$F^o$ and our original 
complex Finsler metric~$F$. The idea is that both measuring the length of 
curves using~$F^o$ and using~$F$ we end up with the same function~$\ell_
\Sigma$; therefore the second variation formula should be the same written in 
terms of~$F$ or in terms of~$F^o$ --- assuming the convexity of the latter, of 
course.
 
The second variation formula for real Finsler metrics has been computed by 
Auslander~[Au1] (see also Chern~[Ch], Bao and Chern~[BC] and~[AP4]), in a 
setting similar to ours and in terms of the so-called horizontal 
flag curvature of the Cartan connection. So comparing the two formulas we get 
an expression for the horizontal flag curvature of the Cartan connection for 
convex K\"ahler Finsler metrics:
 
\newthm Corollary \hfcc: Let $F\colon T^{1,0}M\to\R^+$ be a convex (i.e., with
strongly convex indicatrices) K\"ahler Finsler metric on a complex
manifold~$M$. Then the horizontal flag curvature of the Cartan connection 
associated to~$F^o$ is given by 
$$\eqalign{R_v(H,H)=\Re\Bigl[&\langle\Omega(\chi,\bar H)H,\chi\rangle_v
	-\langle\Omega(H,\smb\chi)H,\chi\rangle_v\cr
	&+\llangle\tau^H(H,\smb\chi),H\rrangle_v-\llangle\tau^H(\chi,\bar H),
	H\rrangle_v\Bigr]\cr}$$
for all $H\in\ca H$.
 
We shall need this result to apply Auslander's version [Au2] of the classical
Cartan-Hadamard theorem. By the way, it turns out that a direct computation of
the Cartan connection (and its curvature) in terms of the Chern-Finsler 
connection (and its curvature) is unexpectedly difficult; see [AP4] for 
details.
 
\smallsect 8. Manifolds with constant holomorphic curvature
 
A very natural problem now is the classification of K\"ahler 
Finsler manifolds of constant holomorphic curvature. In this respect, the 
Finsler situation is much richer than the hermitian one; for instance, 
Lempert's work~[Le] and [AP2] imply that all strongly convex domains of~$\C^n$ 
endowed with the Kobayashi metric are weakly K\"ahler Finsler manifolds with 
constant holomorphic curvature~$-4$.
 
The last theorem of this paper is a step toward this classification; roughly 
speaking, we shall prove that a simply connected K\"ahler Finsler manifold of 
nonpositive constant holomorphic curvature is diffeomorphic to an euclidean 
space. Furthermore, in the case of constant negative holomorphic curvature our 
results show that the Finsler geometry of the manifold is pretty much the same 
of the one of strongly convex domains endowed with the Kobayashi metric.
 
The idea is to apply the Cartan-Hadamard theorem; to do so, we need to 
estimate the curvature terms appearing in the second variation formula.
 
Let $F\colon T^{1,0}M\to\R^+$ be a strongly pseudoconvex Finsler metric on a 
complex manifold~$M$. We say that $F$ has {\sl constant holomorphic curvature} 
$2c\in\R$ if
$$\langle\Omega(\chi,\smb\chi)\chi,\chi\rangle\equiv cG^2,\neweq\eqchc$$
that is iff $K_F\equiv2c$.
The idea is to differentiate \eqchc\ in such a smart way to get all the 
informations we need. 
 
We start with a couple of computational lemmas.
 
\newthm Lemma \ltre: Let $F\colon T^{1,0}M\to\R^+$ be a strongly pseudoconvex 
Finsler metric on a complex manifold~$M$. Then
$$\langle(\nabla_{\bar W}\Omega)(H,\bar K)\chi,\chi\rangle=\langle\tau^H
	\bigl(H,\bar{\theta(K,W)}\bigr),\chi\rangle$$
for all $W\in\ca V$ and $H$,~$K\in\ca H$. In particular,
$$\langle(\nabla_{\bar W}\Omega)(H,\smbar\chi)\chi,\chi\rangle=0$$
for all $W\in\ca V$ and $H\in\ca H$.
 
\pf Since we are interested only in the horizontal part, we may replace 
$\Omega$ by 
$$\Omega^H=\Omega^\alpha_\beta\otimes dz^\beta\otimes\delta_\alpha.$$ 
Since $\nabla_{\bar W}\,dz^\beta=0$ and $\nabla_
{\bar W}\,\delta_\alpha=0$, we have
$$\nabla_{\bar W}\Omega^H=(\nabla_{\bar W}\Omega^\alpha_\beta)
	\otimes dz^\beta\otimes\delta_\alpha.$$
Again, we only need the horizontal part, that is 
$$p^*_H(\nabla_{\bar W}\Omega^\alpha_\beta)=\bar{W}(R^\alpha_{\beta;\mu
	\smb\nu})\,dz^\mu\wedge d\smb z^\nu-R^\alpha_{\beta;\mu\smb\rho}\,\bar{	
	\omega^\rho_\nu(W)}\,dz^\mu\wedge d\smb z^\nu.$$
Recalling \eqR, taking $H$,~$K$,~$L\in\ca H$ we get
$$\eqalign{\langle(&\nabla_{\bar W}\Omega^H(H,\bar K)L,\chi\rangle=
	G_\alpha\bigl[\,\bar{W}(R^\alpha_{\beta;\mu\smb\nu})-R^\alpha_{\beta;\mu
	\smb\rho}\,\bar{\omega^\rho_\nu(W)}\bigr]H^\mu\bar{K^\nu}L^\beta\cr
	&=-G_\alpha\bigl[\dot\de_{\smb\gamma}\delta_{\smb\nu}(\Gamma^\alpha
	_{\beta;\mu})-\delta_{\smb\rho}(\Gamma^\alpha_{\beta;\mu})
	\Gamma^{\smb\rho}
	_{\smb\nu\smb\gamma}+\dot\de_{\smb\gamma}\bigl(\Gamma^\alpha_{\beta
	\sigma}\delta_{\smb\nu}(\Gamma^\sigma_{;\mu})\bigr)-\Gamma^\alpha
	_{\beta\sigma}\delta_{\smb\rho}(\Gamma^\sigma_{;\mu})\Gamma^{\smb\rho}
	_{\smb\nu\smb\gamma}\bigr]H^\mu\bar{K^\nu}L^\beta\bar{W^\gamma}\cr
	&=-G_\alpha\bigl[\delta_{\smb\nu}\dot\de_{\smb\gamma}(\Gamma^\alpha
	_{\beta;\mu})-\Gamma^{\smb\tau}_{\smb\gamma;\smb\nu}\dot\de_{\smb\tau}
	(\Gamma^\alpha_{\beta;\mu})+\dot\de_{\smb\gamma}(\Gamma^\alpha_{\beta
	\sigma})\delta_{\smb\nu}(\Gamma^\sigma_{;\mu})+\delta_{\smb\nu}
	(\Gamma^\sigma_{\smb\gamma;\mu})\Gamma^\alpha_{\beta\sigma}\cr
	&\qquad\qquad-\Gamma
	^\alpha_{\beta\sigma}\Gamma^{\smb\tau}_{\smb\gamma;\smb\nu}\Gamma
	^\sigma_{\smb\tau;\mu}-\delta_{\smb\rho}(\Gamma^\alpha_{\beta;\mu})
	\Gamma^{\smb\rho}_{\smb\nu\smb\gamma}-\Gamma^\alpha_{\beta\sigma}
	\delta_{\smb\rho}(\Gamma^\sigma_{;\mu})\Gamma^{\smb\rho}
	_{\smb\nu\smb\gamma}\bigr]H^\mu\bar{K^\nu}L^\beta\bar{W^\gamma}\cr
\noalign{\noindent(where $\Gamma^\sigma_{\smb\gamma;\mu}=\dot\de_{\smb\gamma}
(\Gamma^\sigma_{;\mu})$ and we used Lemma~\dtre.(ii)),}
	&=-\bigl[\delta_{\smb\nu}\bigl(G_\alpha\dot\de_{\beta}(\Gamma^\alpha
	_{\smb\gamma;\mu})\bigr)-\Gamma^{\smb\tau}_{\smb\gamma;\smb\nu}G_\alpha
	\dot\de_{\beta}(\Gamma^\alpha_{\smb\tau;\mu})+G_{\beta\sigma}\delta_
	{\smb\nu}(\Gamma^\sigma_{\smb\gamma;\mu})-G_{\beta\sigma}\Gamma^{\smb
	\tau}_{\smb\gamma;\smb\nu}\Gamma^\sigma_{\smb\tau;\mu}\cr
	&\qquad-\delta_{\smb\rho}
	(G_\alpha\Gamma^\alpha_{\beta;\mu})\Gamma^{\smb\rho}_{\smb\nu\smb
	\gamma}-
	G_{\beta\sigma}\delta_{\smb\rho}(\Gamma^\sigma_{;\mu})\Gamma^{\smb\rho}
	_{\smb\nu\smb\gamma}\bigr]H^\mu\bar{K^\nu}L^\beta\bar{W^\gamma}\cr
\noalign{\noindent(where we used $G_\alpha\Gamma^\alpha_{\beta\sigma}=
G_{\beta\sigma}$, $\delta_{\smb\nu}(G_\alpha)=0$ and 
$G_\alpha\dot\de_{\smb\gamma}(\Gamma^\alpha_{\beta
\sigma})=G_{\beta\sigma\smb\gamma}-G_{\alpha\smb\gamma}\Gamma^\alpha_{\beta
\sigma}=0$),}
	&=-\bigl[-\delta_{\smb\nu}(G_{\alpha\beta}\Gamma^\alpha
	_{\smb\gamma;\mu})+\Gamma^{\smb\tau}_{\smb\gamma;\smb\nu}G_{\alpha\beta}
	\Gamma^\alpha_{\smb\tau;\mu}+G_{\beta\sigma}\delta_
	{\smb\nu}(\Gamma^\sigma_{\smb\gamma;\mu})-G_{\beta\sigma}\Gamma^{\smb
	\tau}_{\smb\gamma;\smb\nu}\Gamma^\sigma_{\smb\tau;\mu}\cr
	&\qquad-\delta_{\smb\rho}
	(G_\alpha\Gamma^\alpha_{\beta;\mu})\Gamma^{\smb\rho}_{\smb\nu\smb
	\gamma}-
	G_{\beta\sigma}\delta_{\smb\rho}(\Gamma^\sigma_{;\mu})\Gamma^{\smb\rho}
	_{\smb\nu\smb\gamma}\bigr]H^\mu\bar{K^\nu}L^\beta\bar{W^\gamma}\cr
\noalign{\noindent(where we used $G_\alpha\Gamma^\alpha_{\smb\tau;\mu}=0$),}
	&=-\bigl[-\delta_{\smb\nu}(G_{\alpha\beta}\Gamma^\alpha
	_{\smb\gamma;\mu})+G_{\beta\sigma}\delta_{\smb\nu}(\Gamma^\sigma_
	{\smb\gamma;\mu})-\delta_{\smb\rho}(G_\alpha\Gamma^\alpha_{\beta;\mu})
	\Gamma^{\smb\rho}_{\smb\nu\smb\gamma}-
	G_{\beta\sigma}\delta_{\smb\rho}(\Gamma^\sigma_{;\mu})\Gamma^{\smb\rho}
	_{\smb\nu\smb\gamma}\bigr]H^\mu\bar{K^\nu}L^\beta\bar{W^\gamma}.\cr}$$
Hence
$$\eqalign{\langle(&\nabla_{\bar W}\Omega)(H,\bar K)\chi,\chi\rangle\cr
	&=-\bigl[-\delta_{\smb\nu}(G_{\alpha\beta}\Gamma^\alpha
	_{\smb\gamma;\mu})+G_{\beta\sigma}\delta_{\smb\nu}(\Gamma^\sigma_
	{\smb\gamma;\mu})-\delta_{\smb\rho}(G_\alpha\Gamma^\alpha_{\beta;\mu})
	\Gamma^{\smb\rho}_{\smb\nu\smb\gamma}-
	G_{\beta\sigma}\delta_{\smb\rho}(\Gamma^\sigma_{;\mu})\Gamma^{\smb\rho}
	_{\smb\nu\smb\gamma}\bigr]v^\beta H^\mu\bar{K^\nu}\bar{W^\gamma},\cr
	&=\delta_{\smb\rho}(G_\alpha\Gamma^\alpha_{;\mu})\Gamma^{\smb\rho}_
	{\smb\nu\smb\gamma}H^\mu\bar{K^\nu}\bar{W^\gamma},\cr
\noalign{\noindent(where we used \eqdf\ and $v^\beta\Gamma^\alpha_{\beta;\mu}=
\Gamma^\alpha_{;\mu}$),}
	&=G_\alpha\delta_{\smb\rho}(\Gamma^\alpha_{;\mu})\Gamma^{\smb\rho}
	_{\smb\nu\smb\gamma}H^\mu\bar{K^\nu}\,\bar{W^\gamma}=\langle\tau^H
	\bigl(H,\bar{\theta(K,W)}\bigr),\chi\rangle,\cr}$$
because $\theta(K,W)=-\Gamma^\rho_{\nu\gamma}\,K^\nu W^\gamma\delta_\rho$, by 
\eqTau.
 
Finally,
$$\eqalign{\langle(\nabla_{\bar W}\Omega)(H,\smb\chi)\chi,\chi
	\rangle&=G_\alpha\delta_{\smb\rho}(\Gamma^\alpha_{;\mu})
	\Gamma^{\smb\rho}_{\smb\nu\smb\gamma}\bar{v^\nu}H^\mu\bar{W^\gamma}\cr
	&=0,\cr}$$
because $\Gamma^{\smb\rho}_{\smb\nu\smb\gamma}\bar{v^\nu}=0$.\qedn
 
\newthm Lemma \lquat: Let $F\colon T^{1,0}M\to\R^+$ be a strongly pseudoconvex 
Finsler metric on a complex manifold~$M$. Then
$$\langle(\nablab_V\Omega)(H,\bar K)\chi,\chi\rangle=\langle\tau^H
	\bigl(\theta(H,V),\bar K\bigr),\chi\rangle$$
for all $V\in\ca V$ and $H$,~$K\in\ca H$. In particular,
$$\langle(\nabla_V\Omega)(\chi,\bar{K})\chi,\chi\rangle=0$$
for all $V\in\ca V$ and $K\in\ca H$.
 
\pf Again it suffices to consider $\Omega^H=\Omega^\alpha_\beta\otimes dz^\beta
\otimes\delta_\alpha$; so
$$\nabla_V\Omega^H=(\nabla_V\Omega^\alpha_\beta)\otimes dz^\beta
	\otimes\delta_\alpha-\Omega^\alpha_\gamma\otimes\omega^\gamma_\beta(V)
	\,dz^\beta\otimes\delta_\alpha+\Omega^\gamma_\beta\otimes dz^\beta
	\otimes\omega^\alpha_\gamma(V)\delta_\alpha.$$
We are interested only in the horizontal part. Taking $H$,~$K\in\ca H$
we get
$$\eqalign{G_\alpha&(\nablab_V\Omega^\alpha_\beta)(H,\bar{K})=
	G_\alpha\bigl[V(R^\alpha_{\beta;\mu\smb\nu})-R^\alpha_{\beta;\rho
	\smb\nu}\omega^\rho_\nu(V)\bigr]H^\mu\bar{K^\nu}\cr
	&=-G_\alpha\bigl[\dot\de_\lambda\delta_{\smb\nu}(\Gamma^\alpha_{\beta;
	\mu})-\delta_{\smb\nu}(\Gamma^\alpha_{\beta;\rho})\Gamma^\rho_{\mu
	\lambda}+\dot\de_\lambda
	\bigl(\Gamma^\alpha_{\beta\sigma}\delta_{\smb\nu}(\Gamma^\sigma_{;\mu})
	\bigr)-\Gamma^\alpha_{\beta\sigma}\delta_{\smb\nu}(\Gamma^\sigma
	_{;\rho})\Gamma^\rho_{\mu\lambda}\bigr]V^\lambda H^\mu\bar{K^\nu}\cr
	&=-G_\alpha\bigl[\delta_{\smb\nu}\dot\de_\lambda(\Gamma^\alpha_{\beta;
	\mu})-\Gamma^{\smb\tau}_{\lambda;\smb\nu}\dot\de_{\smb\tau}(\Gamma^
	\alpha_{\beta;\mu})-\delta_{\smb\nu}(\Gamma^\alpha_{\beta;\rho})\Gamma^
	\rho_{\mu\lambda}+
	\dot\de_\lambda(\Gamma^\alpha_{\beta\sigma})\delta_{\smb\nu}(\Gamma^
	\sigma_{;\mu})+\Gamma^\alpha_{\beta\sigma}\delta_{\smb\nu}(\Gamma^
	\sigma_{\lambda;\mu})\cr
	&\qquad\qquad-\Gamma^\alpha_{\beta\sigma}\Gamma^\sigma_{\smb\tau;\mu}
	\Gamma^{\smb\tau}_{\lambda;\smb\nu}-\Gamma^\alpha_{\beta\sigma}
	\delta_{\smb\nu}(\Gamma^\sigma_{;\rho})\Gamma^\rho_{\mu\lambda}\bigr]
	V^\lambda H^\mu\bar{K^\nu}\cr
\noalign{\noindent(where we used $[\delta_{\smb\nu},\dot\de_\lambda]=\Gamma^
{\smb\tau}_{\lambda;\smb\nu}\dot\de_{\smb\tau}$),}
	&=-\bigl[\delta_{\smb\nu}\bigl(G_\alpha\dot\de_\beta(\Gamma^\alpha
	_{\lambda;\mu})\bigr)-\Gamma^{\smb\tau}_{\lambda;\smb\nu}G_\alpha
	\dot\de_\beta(\Gamma^\alpha_{\smb\tau;\mu})-\delta_{\smb\nu}(G_\alpha
	\Gamma^\alpha_{\beta;\rho})\Gamma^\rho_{\mu\lambda}+G_\alpha\dot\de
	_\lambda(\Gamma^\alpha_{\beta\sigma})\delta_{\smb\nu}(\Gamma^\sigma
	_{;\mu})\cr
	&\qquad+G_{\beta\sigma}\delta
	_{\smb\nu}(\Gamma^\sigma_{\lambda;\mu})-G_{\beta\sigma}\Gamma^\sigma
	_{\smb\tau;\mu}\Gamma^{\smb\tau}_{\lambda;\smb\nu}-G_{\beta\sigma}
	\delta_{\smb\nu}(\Gamma^\sigma_{;\rho})\Gamma^\rho_{\mu\lambda}\bigr]
	V^\lambda H^\mu\bar{K^\nu}\cr
\noalign{\noindent(where we used $\delta_{\smb\nu}(G_\alpha)=0$ and 
$G_\alpha\Gamma^\alpha_{\beta\sigma}=G_{\beta\sigma}$),}
	&=-\bigl[-\delta_{\smb\nu}(G_{\alpha\beta}\Gamma^\alpha_{\lambda;\mu})
	+\delta_{\smb\nu}\bigl(\dot\de_\beta\delta_\mu(G_\lambda)\bigr)+\Gamma
	^{\smb\tau}_{\lambda;\smb\nu}G_{\alpha\beta}\Gamma^\alpha_{\smb\tau;\mu}
	-\delta_{\smb\nu}(G_\alpha\Gamma^\alpha_{\beta;\rho})\Gamma^\rho_{\mu
	\lambda}\cr
	&\qquad+
	G_\alpha\dot\de_\lambda(\Gamma^\alpha_{\beta\sigma})\delta_{\smb\nu}
	(\Gamma^\sigma_{;\mu})+G_{\beta\sigma}\delta
	_{\smb\nu}(\Gamma^\sigma_{\lambda;\mu})-G_{\beta\sigma}\Gamma^\sigma
	_{\smb\tau;\mu}\Gamma^{\smb\tau}_{\lambda;\smb\nu}-G_{\beta\sigma}
	\delta_{\smb\nu}(\Gamma^\sigma_{;\rho})\Gamma^\rho_{\mu\lambda}\bigr] 
	V^\lambda H^\mu\bar{K^\nu}\cr
\noalign{\noindent(where we used $G_\alpha\Gamma^\alpha_{\smb\tau;\mu}=0$ and 
$G_\alpha\Gamma^\alpha_{\lambda;\mu}=\delta_\mu(G_\lambda)$),}
	&=-\bigl[-\delta_{\smb\nu}(G_{\alpha\beta}\Gamma^\alpha_{\lambda;\mu})
	+\delta_{\smb\nu}\bigl(\dot\de_\beta\delta_\mu(G_\lambda)\bigr)
	-\delta_{\smb\nu}(G_\alpha\Gamma^\alpha_{\beta;\rho})\Gamma^\rho_{\mu
	\lambda}+G_\alpha\dot\de_\lambda(\Gamma^\alpha_{\beta\sigma})\delta
 	_{\smb\nu}(\Gamma^\sigma_{;\mu})\cr
	&\qquad+G_{\beta\sigma}\delta
	_{\smb\nu}(\Gamma^\sigma_{\lambda;\mu})-G_{\beta\sigma}
	\delta_{\smb\nu}(\Gamma^\sigma_{;\rho})\Gamma^\rho_{\mu\lambda}\bigr]
	V^\lambda H^\mu\bar{K^\nu}.\cr}$$
Furthermore,
$$\eqalign{G_\alpha\Omega^\alpha_\gamma(H,\bar{K})\,\omega^\gamma_\beta(V)&=
	-G_\alpha\bigl[\delta_{\smb\nu}(\Gamma^\alpha_{\gamma;\mu})+\Gamma
	^\alpha_{\gamma\sigma}\delta_{\smb\nu}(\Gamma^\sigma_{;\mu})\bigr]
	\Gamma^\gamma_{\beta\lambda}V^\lambda H^\mu\bar{K^\nu}\cr
	&=-\bigl[\delta_{\smb\nu}\bigl(\delta_\mu(G_\gamma)\bigr)+G_{\gamma
	\sigma}\delta_{\smb\nu}(\Gamma^\sigma_{;\mu})\bigr]\Gamma^\gamma
	_{\beta\lambda}V^\lambda H^\mu\bar{K^\nu};\cr}$$
$$\eqalign{G_\alpha\,\omega^\alpha_\gamma(V)\Omega^\gamma_\beta(H,\bar{K})&=
	-G_\alpha\Gamma^\alpha_{\gamma\lambda}\bigl[\delta_{\smb\nu}(\Gamma
	^\gamma_{\beta;\mu})+\Gamma^\gamma_{\beta\sigma}\delta_{\smb\nu}
	(\Gamma^\sigma_{;\mu})\bigr]V^\lambda H^\mu\bar{K^\nu}\cr
	&=-G_{\gamma\nu}\bigl[\delta_{\smb\nu}(\Gamma
	^\gamma_{\beta;\mu})+\Gamma^\gamma_{\beta\sigma}\delta_{\smb\nu}
	(\Gamma^\sigma_{;\mu})\bigr]V^\lambda H^\mu\bar{K^\nu}.\cr}$$
Summing up we find
$$\eqalign{\langle&(\nablab_V\Omega)(H,\bar{K})\chi,\chi\rangle\cr
	&=-\bigl[-\delta_{\smb\nu}(G_{\alpha\beta}\Gamma^\alpha_{\lambda;\mu})
	+\delta_{\smb\nu}\bigl(\dot\de_\beta\delta_\mu(G_\lambda)\bigr)
	-\delta_{\smb\nu}(G_\alpha\Gamma^\alpha_{\beta;\rho})\Gamma^\rho_{\mu
	\lambda}+G_\alpha\dot\de_\lambda(\Gamma^\alpha_{\beta\sigma})
	\delta_{\smb\nu}(\Gamma^\sigma_{;\mu})\cr
	&\qquad+G_{\beta\sigma}\delta_{\smb\nu}(\Gamma^\sigma_{\lambda;\mu})
	-G_{\beta\sigma}\delta_{\smb\nu}(\Gamma^\sigma_{;\rho})\Gamma^\rho
	_{\mu\lambda}+\delta_{\smb\nu}\bigl(\delta_\mu(G_\gamma)\bigr)
	\Gamma^\gamma_{\beta\lambda}
	+G_{\gamma\sigma}\delta_{\smb\nu}(\Gamma^\sigma_{;\mu})\Gamma^\gamma
	_{\beta\lambda}+G_{\gamma\nu}\delta_{\smb\nu}(\Gamma
	^\gamma_{\beta;\mu})\cr
	&\qquad+G_{\gamma\nu}\Gamma^\gamma_{\beta\sigma}\delta
	_{\smb\nu}(\Gamma^\sigma_{;\mu})\bigr]v^\beta V^\lambda H^\mu
	\bar{K^\nu}\cr
	&=-\bigl[\delta_{\smb\nu}\bigl(v^\beta\dot\de_\beta\delta_\mu(G_\lambda)
	\bigr)-\delta_{\smb\nu}(G_\alpha v^\beta\Gamma^\alpha_{\beta;\rho})
	\Gamma^\rho_{\mu\lambda}+G_\alpha v^\beta\dot\de_\lambda(\Gamma^\alpha_
	{\beta\sigma})\delta_{\smb\nu}(\Gamma^\sigma_{;\mu})\cr
	&\qquad+G_{\gamma\lambda}\delta_{\smb\nu}
	(\Gamma^\gamma_{;\mu})\bigr]V^\lambda H^\mu\bar{K^\nu}\cr
	&=-\bigl[-\delta_{\smb\nu}(G_{;\rho})\Gamma^\rho_{\mu\lambda}
	-G_{\lambda\sigma}
	\delta_{\smb\nu}(\Gamma^\sigma_{;\mu})+G_{\lambda\gamma}\delta_{\smb\nu}
	(\Gamma^\gamma_{;\mu})\bigr]V^\lambda H^\mu\bar{K^\nu}\cr
\noalign{\noindent(where we used $v^\beta\dot\de_\beta\delta_\mu(G_\lambda)=
v^\beta G_{\lambda\beta;\mu}-v^\beta\Gamma^\sigma_{\beta;\mu}G_{\lambda\sigma}-
v^\beta\Gamma^\sigma_{;\mu}G_{\lambda\sigma\beta}=0$),}
	&=\delta_{\smb\nu}(G_{;\rho})\Gamma^\rho_{\mu\lambda}V^\lambda H^\mu
	\bar{K^\nu}=
	\delta_{\smb\nu}(G_\sigma\Gamma^\sigma_{;\rho})\Gamma^\rho_{\mu\lambda}
	V^\lambda H^\mu\bar{K^\nu}\cr
	&=G_\sigma\delta_{\smb\nu}(\Gamma^\sigma_{;\rho})\Gamma^\rho_{\mu
	\lambda}
	V^\lambda H^\mu\bar{K^\nu}=\langle\tau^H\bigl(\theta(H,V),\bar{K}\bigr),
	\chi\rangle.\cr}$$
The final assertion follows from $\theta(\chi,V)=0$.\qedn
 
In the following computations we shall need some symmetries of the curvature 
operator, summarized in 
 
\newthm Lemma \ldue: Let $F\colon T^{1,0}M\to\R^+$ be a strongly pseudoconvex 
Finsler metric on a complex manifold~$M$. Then
{\smallskip
\itm{(i)}$\langle\Omega(H,\smbar\chi)\chi,\chi\rangle=\langle\Omega(\chi,\smbar
\chi)H,\chi\rangle$ for all $H\in\ca H$ iff
$$\langle\smb\de_H\theta(H,\chi,\smbar\chi),\chi\rangle=0\neweq\eqsuno$$ 
for all $H\in\ca H$;
\itm{(ii)}$\langle\Omega(H,\bar K)\chi,\chi\rangle=\langle\Omega(\chi,\bar K)H,
\chi\rangle$ for all $H$,~$K\in\ca H$ iff
$$\langle\smb\de_H \theta(H,\chi,\bar K),\chi\rangle=0\neweq\eqsdue$$
for all $H$,~$K\in\ca H$.}
 
\pf It follows immediately from \eqtth\ and Proposition~\dB.\qedn
 
Now we can start. The first step is:
 
\newthm Proposition \lcin: Let $F\colon T^{1,0}M\to\R^+$ be a strongly 
pseudoconvex Finsler metric on a complex manifold~$M$, with constant 
holomorphic curvature~$2c\in\R$. Then
$$\langle\smb\de_H\theta(H,\chi,\smbar\chi),\chi\rangle=0\neweq\eqHzero$$
for all $H\in\ca H$ iff
$$\tau^H(\chi,\smbar\chi)=cG\chi.\neweq\eqstar$$
Furthermore, they both imply
$$\langle\Omega(\chi,\bar{K})\chi,\chi\rangle=cG\langle\chi,K\rangle
	\neweq\eqCuno$$
for all $K\in\ca H$.
 
\pf Take $W\in\ca V$ and let $K=\Theta(W)\in\ca H$; note that $\nabla_{\bar 
W}\raise2pt\hbox{$\chi$}=0$ and $\nablab_W\raise2pt\hbox{$\chi$}=\Theta(W)=K$. 
Then
$$\eqalign{\bar{W}(cG^2)&=2cG\langle\chi,K\rangle;\cr
	\bar{W}\langle\Omega(\chi,\smbar\chi)\chi,\chi\rangle&=\langle
	(\nabla_{\bar W}\Omega)(\chi,\smbar\chi)\chi,\chi\rangle+
	\langle\Omega(\chi,\bar{K})\chi,\chi\rangle+\langle\Omega(\chi,
	\smbar\chi)\chi,K\rangle\cr
	&=\langle\Omega(\chi,\bar{K})\chi,\chi\rangle+\langle\tau^H(\chi,
	\smbar\chi),K\rangle,\cr}\neweq\eqaddi$$
where we used Lemmas~\ltre\ and~\luno. Since $F$ has constant holomorphic 
curvature~$2c$, we have
$$\langle\Omega(\chi,\smb\chi)\chi,\chi\rangle=cG^2$$
and hence \eqaddi\ yields
$$\langle\Omega(\chi,\bar{K})\chi,\chi\rangle=2cG\langle\chi,K\rangle
	-\langle\tau^H(\chi,\smbar\chi),K\rangle.\neweq\eqstarc$$
Subtracting $\langle\tau^H(\chi,\smbar\chi),K\rangle=\langle\Omega(\chi,\smb
\chi)\chi,K\rangle$ to both sides, we find that \eqstar\ holds if and only if
$$\langle\Omega^H(\chi,\bar{K})\chi,\chi\rangle=\langle\Omega^H(\chi,
	\smbar\chi)\chi,K\rangle$$
for all $K\in\ca H$, that is, recalling \eqtcb, iff
$$\langle\Omega(K,\smbar\chi)\chi,\chi\rangle=\langle\Omega(\chi,\smbar\chi)
	K,\chi\rangle,$$
and thus, by Lemma~\ldue, iff \eqHzero\ holds.
 
Finally, if \eqstar\ holds, \eqstarc\ yields \eqCuno.\qedn
 
The second step requires \eqsdue:
 
\newthm Proposition \tuuno: Let $F\colon T^{1,0}M\to\R^+$ be a strongly 
pseudoconvex Finsler metric on a complex manifold~$M$ with constant 
holomorphic curvature~$2c\in\R$. Assume that $\eqsuno$ holds. Then
$$\langle\Omega(H,\bar K)\chi,\chi\rangle+\langle\Omega(\chi,\bar K)H,\chi
	\rangle=c\left\{\langle H,\chi\rangle\langle\chi,K\rangle+\langle\chi,
	\chi\rangle\langle H,K\rangle\right\},\neweq\eqCdue$$
for all $H$,~$K\in\ca H$. In particular, if $\eqsdue$ holds then
$$\langle\Omega(\chi,\bar K)H,\chi\rangle={c\over2}\left\{\langle H,\chi\rangle
	\langle\chi,K\rangle+\langle\chi,\chi\rangle\langle H,K\rangle\right\}
	\neweq\eqHbis$$
for all $H$,~$K\in\ca H$.
 
\pf Take $V$,~$W\in\ca V$ such that $\Theta(V)=H$ and $\Theta(W)=K$ and extend 
them in any way to sections of $\ca V$ (and thus extend $H$ and $K$ as 
sections of $\ca H$ via $\Theta$).
We have
$$V\bigl(cG\langle\chi,K\rangle)=c\bigl[\langle H,\chi\rangle\langle\chi,K
	\rangle+G\langle H,K\rangle+G\langle\chi,\nabla_{\bar V}K\rangle\bigr]
	,$$
and
$$\eqalign{V\langle&\Omega(\chi,\bar K)\chi,\chi\rangle\cr
	&=\langle(\nablab_V\Omega)(\chi,\bar K)\chi,\chi\rangle+\langle
	\Omega(H,\bar K)\chi,\chi\rangle+\langle\Omega(\chi,\bar{
	\nablab_VK})\chi,\chi\rangle+\langle\Omega(\chi,\bar K)H,\chi\rangle\cr
	&=\langle\Omega(H,\bar K)\chi,\chi\rangle+\langle\Omega(\chi,
	\bar{\nabla_{\bar V}K})\chi,\chi\rangle+\langle\Omega(\chi,
	\bar K)H,\chi\rangle,\cr}$$
thanks to Lemma~\lquat. Since \eqsuno\ holds, we can use Lemma~\lcin\ (that 
is, \eqCuno\ applied both to $K$ and to $\nabla_{\bar V}K$) to get
exactly \eqCdue.
 
Finally, \eqHbis\ follows from Lemma~\ldue.\qedn
 
So we have obtained one of the hermitian product terms. This immediately 
yields one of the symmetric product terms:
 
\newthm Proposition \tstep: Let $F\colon T^{1,0}M\to\R^+$ be a strongly 
pseudoconvex Finsler metric on a complex manifold~$M$ with constant 
holomorphic curvature~$2c\in\R$. Assume that $\eqsdue$ holds. Then
$$\tau^H(K,\smb\chi)={c\over2}\left\{\langle K,\chi\rangle\chi+\langle\chi,
	\chi\rangle K\right\}\neweq\eqtsu$$
for all $K\in\ca H$. In particular,
$$\llangle H,\tau^H(K,\smb\chi)\rrangle={c\over2}\langle\chi,\chi\rangle
	\llangle H,K\rrangle\neweq\eqtsd$$
for all $H$,~$K\in\ca H$.
 
\pf We get
$$\langle H,\tau^H(K,\smb\chi)\rangle=\langle H,\Omega(K,\smb\chi)\chi\rangle
	=\langle\Omega(\chi,\bar K)H,\chi\rangle$$
for all $H$,~$K\in\ca H$, thanks to Lemma~\luno\ and \eqtcb. Then \eqHbis\ 
yields \eqtsu, and \eqtsd\ follows immediately.\qedn
 
For the other symmetric product term we need the weak K\"ahler condition:
 
\newthm Proposition \fstep: Let $F\colon T^{1,0}M\to\R^+$ be a weakly 
K\"ahler-Finsler metric on a complex manifold~$M$ such that $\eqsdue$ holds. 
Then
$$\llangle H,\tau^H(\chi,\bar K)\rrangle=0$$
for all $H$,~$K\in\ca H$.
 
\pf The weak K\"ahler condition $\langle \theta(H,\chi),\chi\rangle=0$ for all 
$H\in\ca H$ implies
$$\forevery{H,\,K\in\ca H}\langle(\nabla_{\bar K}\theta)(H,\chi),\chi\rangle=0
	\neweq\eqbKH$$
because $\nabla_{\bar K}\raise2pt\hbox{$\chi$}=0=\nablab_K\raise2pt
\hbox{$\chi$}$. Now, writing $\theta=\theta^
\alpha\otimes\delta_\alpha$, we have $\nabla_{\bar K}\theta=(\nabla_{\bar K}
\theta^\alpha)\otimes\delta_\alpha$ and
$$\nabla_{\bar K}\theta^\alpha=\bar{K^\tau}\delta_{\smb\tau}(\Gamma^\alpha_
	{\nu;\mu})\,dz^\mu\wedge dz^\nu+\bar{K^\tau}\delta_{\smb\tau}(\Gamma^
	\alpha_{\nu\gamma})\,\psi^\gamma\wedge dz^\nu.$$
Therefore \eqbKH\ implies
$$G_\alpha[\delta_{\smb\tau}(\Gamma^\alpha_{\nu;\mu})-\delta_{\smb\tau}
	(\Gamma^\alpha_{\mu;\nu})]H^\mu\bar{K^\tau}v^\nu=0
	\neweq\eqbKHb$$
for all $H$,~$K\in\ca H$.
 
Writing the curvature in local coordinates we find
$$\eqalign{\langle\Omega(\chi,\bar K)H,\chi\rangle&=-G_\alpha[\delta_{\smb\tau}
	(\Gamma^\alpha_{\mu;\nu})+\Gamma^\alpha_{\mu\sigma}\delta_{\smb\tau}
	(\Gamma^\sigma_{;\nu})]H^\mu\bar{K^\tau}v^\nu,\cr
	\langle\Omega(H,\bar K)\chi,\chi\rangle&=-G_\alpha[\delta_{\smb\tau}
	(\Gamma^\alpha_{\nu;\mu})+\Gamma^\alpha_{\nu\sigma}\delta_{\smb\tau}
	(\Gamma^\sigma_{;\mu})]H^\mu\bar{K^\tau}v^\nu.\cr}$$
So \eqbKHb\ yields
$$\eqalign{\langle\Omega(\chi,\bar K)H,\chi\rangle-\langle\Omega(H,\bar K)\chi,
	\chi\rangle&=-G_\alpha\Gamma^\alpha_{\mu\sigma}\delta_{\smb\tau}
	(\Gamma^\sigma_{;\nu})H^\mu\bar{K^\tau}v^\nu\cr
	&=\llangle H,\tau^H(\chi,\bar K)\rrangle,\cr}$$
and the assertion follows from \eqsdue.\qedn
 
We are left with the last term:
 
\newthm Proposition \vstep: Let $F\colon T^{1,0}M\to\R^+$ be a strongly 
pseudoconvex Finsler metric on a complex manifold~$M$ with constant 
holomorphic curvature~$2c\in\R$. Assume that $\eqsuno$ holds. 
Then
$$\langle\Omega(H,\smb\chi)K,\chi\rangle=c\left\{\langle H,\chi\rangle
	\langle K,\chi\rangle+\langle\chi,\chi\rangle\llangle H,K\rrangle
	\right\}$$
for all $H$,~$K\in\ca H$.
 
\pf First of all, we have
$$\langle\Omega(H,\smb\chi)\chi,\chi\rangle=\langle\chi,\Omega(\chi,\bar H)
	\chi\rangle=\bar{\langle\Omega(\chi,\bar H)\chi,\chi\rangle}
	=cG\langle H,\chi,\rangle,\neweq\eqbboh$$
by \eqtcb\ and \eqCuno. Now take $W\in\ca V$ such that $\Theta(W)=K$; then
$$\eqalign{W(cG\langle H,\chi\rangle)&=c\left\{\langle K,\chi\rangle\langle
	H,\chi\rangle+G\langle\nablab_WH,\chi\rangle\right\},\cr
	W\langle\Omega(H,\smb\chi)\chi,\chi\rangle&=\langle(\nablab_W\Omega)
	(H,\smb\chi)\chi,\chi\rangle+\langle\Omega(\nablab_WH,\smb\chi)\chi,
	\chi\rangle+\langle\Omega(H,\smb\chi)K,\chi\rangle,\cr}$$
and so \eqbboh\ yields
$$\langle\Omega(H,\smb\chi)K,\chi\rangle=c\langle K,\chi\rangle\langle H,\chi
	\rangle-\langle(\nablab_W\Omega)(H,\smb\chi)\chi,\chi\rangle.$$
Now Lemma~\lquat\ gives
$$\eqalign{\langle(\nablab_W\Omega)(H,\smb\chi)\chi,\chi\rangle&=
	\langle\tau^H\bigl(\theta(H,W),\smb\chi\bigr),\chi\rangle=
	\langle\Omega\bigl(\theta(H,W),\smb\chi)\chi,\chi\rangle\cr
	&=cG\langle \theta(H,W),\chi\rangle,\cr}$$
again by \eqbboh. But
$$\langle \theta(H,W),\chi\rangle=-G_\alpha\Gamma^\alpha_{\nu\beta}K^\beta
	H^\nu=-\llangle H,K\rrangle,$$
and we are done.\qedn
 
We can finally collect all our computations in
 
\newthm Corollary \fine: Let $F\colon T^{1,0}M\to\R^+$ be a weakly 
K\"ahler-Finsler metric on a complex manifold~$M$. Assume $F$ has constant 
holomorphic curvature~$2c\in\R$ and that $\eqsdue$ holds. Then
$$\displaylines{\qquad\Re\Bigl[\langle\Omega(\chi,\bar K)H,\chi\rangle-
	\langle\Omega(H,\smb\chi)K,\chi\rangle+\llangle H,\tau^H(K,\smb\chi)
	\rrangle-\llangle H,\tau^H(\chi,\bar K)\rrangle\Bigr]\hfill\cr
	\hfill={c\over2}\Re\Bigl[G\bigl\{\langle H,K\rangle-\llangle H,K
	\rrangle\bigr\}+\langle H,\chi\rangle\bigl\{\langle\chi,K\rangle-
	2\langle K,\chi\rangle\bigr\}\Bigr]\qquad\cr}$$
for all $H$,~$K\in\ca H$. 
 
\pf It follows from Propositions~\tuuno, \tstep, \fstep, \vstep\ and 
Corollary~\hfcc.\qedn
 
We are then able to prove the announced
 
\newthm Theorem \final: Let $F\colon T^{1,0}M\to\R^+$ be a complete Finsler 
metric on a simply connected complex manifold~$M$. Assume that:
{\smallskip
\itm{(i)} $F$ is K\"ahler;
\itm{(ii)} $F$ has nonpositive constant holomorphic curvature $2c\le0$;
\itm{(iii)} $\langle\smb\de_H\theta(H,\chi,\bar K),\chi\rangle=0$ for 
all~$H$,~$K\in\ca H$;
\itm{(iv)} the indicatrices of $F$ are strongly convex.
\smallskip
\noindent Then $\exp_p\colon T^{1,0}_p\to M$ is a homeomorphism, and a smooth 
diffeomorphism outside the origin, for any~$p\in M$. Furthermore, $M$ is 
foliated by isometric totally geodesic holomorphic embeddings of the unit 
disk~$\Delta$ endowed with a suitable multiple of the Poincar\'e metric if $c<
0$, or by isometric totally geodesic holomorphic embeddings of~$\C$ endowed 
with the euclidean metric if~$c=0$. In particular, if $2c=-4$ then $F$ is the 
Kobayashi metric of~$M$, and if $c=0$ then the Kobayashi metric of~$M$ 
vanishes identically.}
 
\pf Let $F^o\colon T_\R M\to\R^+$ be the real Finsler metric associated to~$F$
as at the beginning of section~6. Then Corollary~\hfcc\ and Corollary~\fine\
show that the horizontal flag curvature of~$F^o$ is given by 
$$R(H,H)={c\over 2}\Re\Bigl\{G\bigl[\langle H,H\rangle-\llangle H,H\rrangle
	\bigr]+\langle H,\chi\rangle\bigl[\langle\chi,H\rangle-2\langle H,
	\chi\rangle\bigr]\Bigr\}.$$
In particular, if $H=\chi$ we get
$$R(\chi,\chi)=0,\neweq\eqqu$$
and if $\langle H,\chi\rangle=0$ we get
$$R(H,H)={cG\over2}\Re\bigl[\langle H,H\rangle-\llangle H,H\rrangle\bigr]
	={cG\over2}\Re\bigl[\langle iH,iH\rangle+\llangle iH,iH\rrangle\bigr]
	.$$
Now, in local coordinates the quadratic form
$$H\mapsto\Re\bigl[\langle H,H\rangle+\llangle H,H\rrangle\bigr]$$
is represented by the Hessian of~$G$; by (iv), it is positive definite. So
$\langle H,\chi\rangle$ implies
$$R(H,H)\le0.\neweq\eqqdu$$
Now, take $K\in\ca H$ and write $K=\zeta\chi+H$, with $\langle H,\chi\rangle=
0$. Then Corollary~\fine, \eqqu\ and \eqqdu\ yield
$$R(K,K)=R(H,H)\le 0.$$
In conclusion, the horizontal flag curvature is negative semi-definite, and 
the first assertion follows from Auslander's version of the Cartan-Hadamard
theorem~[Au2]. Finally, the last assertion has been proved under  
weaker assumptions in~[AP2,~3].\qedn
 
We remark that, contrarily to what happens in the hermitian 
case, condition (iii) does not seem to be a consequence of the K\"ahler 
condition. For instance, the proof of Proposition~\fstep\ shows that if $F$ is 
weakly K\"ahler (but even K\"ahler does not help much) then condition~(iii) 
holds iff
$$\forevery{H,\,K\in\ca H}\llangle H,\tau^H(\chi,\bar K)\rrangle=0.$$
 
\setref{ZZZ}\nobreak
\beginsection References
 
\book A M. Abate: Iteration theory of holomorphic maps on taut manifolds! 
Me\-di\-ter\-ra\-nean Press, Cosenza, 1989
 
\art AP1 M. Abate, G. Patrizio: Uniqueness of complex geodesics and 
characterization of circular domains! Man. Math.! 74 1992 277-297
 
\book AP2 M. Abate, G. Patrizio: {\sl Holomorphic curvature of Finsler metrics 
and complex geodesics}! To appear in J. Geom. Anal. (1993)
 
\book AP3 M. Abate, G. Patrizio: {\sl Complex geodesics and Finsler metrics}!
To appear in Proceedings of the Semester in Complex Analysis, Banach Center,
Warsaw, 1992
 
\book AP4 M. Abate, G. Patrizio: {\sl A global approach to Finsler metrics}!
In preparation (1993)
 
\art Au1 L. Auslander: On the use of forms in the variational calculations! 
Pac. J. Math.! 5 1955 853-859
 
\art Au2 L. Auslander: On curvature in Finsler geometry! Trans. Am. Math. 
Soc.! 79 1955 378-388
 
\art BC D. Bao, S.S. Chern: On a notable connection in Finsler geometry! 
Houston J. Math.! 19 1993 138-180
 
\book B A. Bejancu: Finsler geometry and applications! Ellis Horwood Limited,
Chichester, 1990
 
\book C E. Cartan: Les espaces de Finsler! Hermann, Paris, 1934
 
\art Ch S.S. Chern: On Finsler geometry! C.R. Acad. Sc. Paris! 314 1992 757-761
 
\art F J.J. Faran: Hermitian Finsler metrics and the Kobayashi metric! J. 
Diff. Geom.! 31 1990 601-625
 
\art Fu M. Fukui: Complex Finsler manifolds! J. Math. Kyoto Univ.! 29 1989 
609-624
 
\book JP M. Jarnicki, P. Pflug: Invariant distances and metrics in complex 
analysis! de Gruyter, Berlin, 1993
 
\art K1 S. Kobayashi: Invariant distances on complex manifolds and holomorphic 
mappings! J. Math. Soc. Japan! 19 1967 460-480
 
\book K2 S. Kobayashi: Hyperbolic manifolds and holomorphic mappings! Dekker, 
New York, 1970
 
\art K3 S. Kobayashi: Negative vector bundles and complex Finsler structures! 
Nagoya Math. J.! 57 1975 153-166
 
\art K4 S. Kobayashi: Intrinsic distances, measures and geometric function 
theory! Bull. Am. Math. Soc.! 82 1976 357-416
 
\book KN S. Kobayashi, K. Nomizu: Foundations of differential geometry, 
vol.~I\negthinspace I! Wiley, New York, 1969
 
\book L S. Lang: An introduction to complex hyperbolic spaces! Springer, New 
York, 1987
 
\art Le L. Lempert: La m\'etrique de Kobayashi et la repr\'esentation des 
domaines sur la boule! Bull. Soc. Math. France! 109 1981 427-474
 
\book M M. Matsumoto: Foundations of Finsler geometry and special Finsler 
spaces! Kaiseisha Press, $\bar{\hbox{O}}$tsu Japan, 1966
 
\art P1 M.Y. Pang: Finsler metrics with the properties of the Kobayashi 
metric on convex domains! Publications Math\'ematiques! 36 1992 131-155
 
\book P2 M.Y. Pang: {\sl Smoothness of the Kobayashi metric of non-convex 
domains}! Preprint (1993)
 
\book R H.L. Royden: {\rm Complex Finsler metrics}! In {\bf Contemporary 
Mathematics. Proceedings of Summer Research Conference,} American Mathematical 
Society, Providence, 1984, pp.~119--124
 
\book Ru1 H. Rund: The differential geometry of Finsler spaces! Springer, 
Berlin, 1959
 
\art Ru2 H. Rund: Generalized metrics on complex manifolds! Math. Nach.! 34 
1967 55-77
 
\art S M. Suzuki: The intrinsic metrics on the domains in $\C^n$! Math. Rep. 
Toyama Univ.! 6 1983 143-177
 
\art W B. Wong: On the holomorphic sectional curvature of some intrinsic 
metrics! Proc. Am. Math. Soc.! 65 1977 57-61
 
\art Wu H. Wu: A remark on holomorphic sectional curvature! Indiana Math. J.! 
22 1973 1103-1108
 
\vfill	
\bye